\documentclass[a4paper,11pt]{article}

\RequirePackage[OT1]{fontenc} \RequirePackage{amsthm,amsmath,natbib}
\RequirePackage[colorlinks]{hyperref} \RequirePackage{hypernat}
\usepackage{amssymb}
\usepackage{graphicx}
\usepackage{latexsym}
\usepackage{amsfonts}
\usepackage[english]{babel}
\usepackage{xy}
\usepackage{fancyhdr}

\usepackage{amssymb}
\usepackage{natbib}

\usepackage{latexsym,enumerate}
\usepackage{amsmath,amsthm,amstext,amssymb}

\newcommand{\Z}{\ensuremath{\mathbb{Z}}}

\newtheorem{theorem}{Theorem}

\newtheorem{proposition}{Proposition}
\newtheorem{lemma}{Lemma}
\newtheorem{definition}{Definition}

\newdimen\AAdi%
\newbox\AAbo%
\def\AArm{\fam0 \mathrm}%
\def\AAk#1#2{\setbox\AAbo=\hbox{#2}\AAdi=\wd\AAbo\kern#1\AAdi{}}%
\def\AAr#1#2#3{\setbox\AAbo=\hbox{#2}\AAdi=\ht\AAbo\raise#1\AAdi\hbox{#3}}%
\def\BBn{{\AArm I\!N}}%
\def\BBr{{\AArm I\!R}}%
\def\BBz{{\AArm Z\!\!Z}}%
\def \P{\hbox{\it I\!P}}
\def \E{\hbox{\it I\hskip -2pt E}}
\def \V{\mbox{\rm Var}}
\def \I{\hbox{\rm 1\hskip -3pt I}}
\def \I{\hbox{\rm 1\hskip -3pt I}}

\begin{document}

\begin{center}
{\bf {\sc \Large {A test of goodness-of-fit for the copula
densities}}}

\vspace{1cm}

Ghislaine GAYRAUD $^*$ and Karine TRIBOULEY $^\flat$

\vspace{0.5cm}

$^*$ CREST and LMRS-UMR 6085\\ Universit\'e de Rouen\\
 Avenue de l'Universit\'e, BP.12\\
  76801 Saint-\'Etienne-du-Rouvray, FRANCE\\
\smallskip
\textsf{Ghislaine.Gayraud@univ-rouen.fr}\\
\bigskip

$^\flat$ LPMA and ModalX\\Universite Paris 10\\
Batiment G\\
 200 avenue de la R\'epublique\\
 92 001 Nanterre Cedex,  FRANCE\\
\smallskip
\textsf{karine.tribouley@u-paris10.fr}

\vspace{0.5cm}

\end{center}

\vspace{0.5cm}

\begin{abstract}
\noindent {\rm \hspace{0.2cm} We consider the problem of testing
hypotheses on the copula density from $n$ bi-dimensional
observations. We wish to test the null hypothesis characterized by a
parametric class against a composite nonparametric alternative. Each
density under the alternative is separated in the $L_2$-norm from
any density lying in the null hypothesis. The copula densities under
consideration are supposed to belong to a range of Besov balls.
According to the minimax approach, the testing problem is solved in
an adaptive framework: it leads to a $\log\log$ term loss in the
minimax rate of testing in comparison with the non-adaptive case. A
smoothness-free test statistic that achieves the minimax rate is
proposed.  The lower bound is also proved. Besides, the empirical
performance of the test procedure is demonstrated with both
simulated and real data.\\

\vspace{0.5cm}



}
\end{abstract}

\newpage

\section{Introduction}
Copulas became a very popular and attractive tool in the recent
literature for modeling multivariate observations. The nice
feature of copulas is that they capture the structure dependence among the components of a multivariate observation without
requiring the study of the univariate margins. More precisely, Sklar's Theorem ensures that any $d-$varied
distribution function $H$ may be expressed as
$$
H(x^1,\ldots,x^d)=C\left(F^1(x^1),\ldots,F^d(x^d)\right),
$$
where the $F^p$'s are the margins and $C$ is called the copula
function. \citep{Sklar:1959} states the existence and the uniqueness
of $C$ as soon as the random variables with joint law $H$ are
continuous.

Modeling the dependence is a great challenge in statistics,
specially in finance or assurance where (for instance) the
identification of the dependence structure between assets is
essential. Many authors proposed parametrical families of copulas
$\{C_\lambda,\lambda\in \Lambda\}$, each of them being available to
capture different dependence behavior. The elliptic family contains
the Gaussian copulas and the Student copula which are often used in
finance. For insurance purposes, heavy tails are needed and copulas
coming from the archimedian family are used. Among others, the more
common are the Gumbel copula, the Clayton copula or the Frank
copula. In view to illustrate the different behaviours of the tails
of several copula densities, some graphs corresponding to the
models cited above are presented below. The parameters are chosen such a way that the
associated Kendall's tau (i.e. the indicator of
concordance/discordance) is identical in all illustrations.

\begin{figure}[h!]
\begin{center}
\includegraphics[width=4cm,height=4cm]{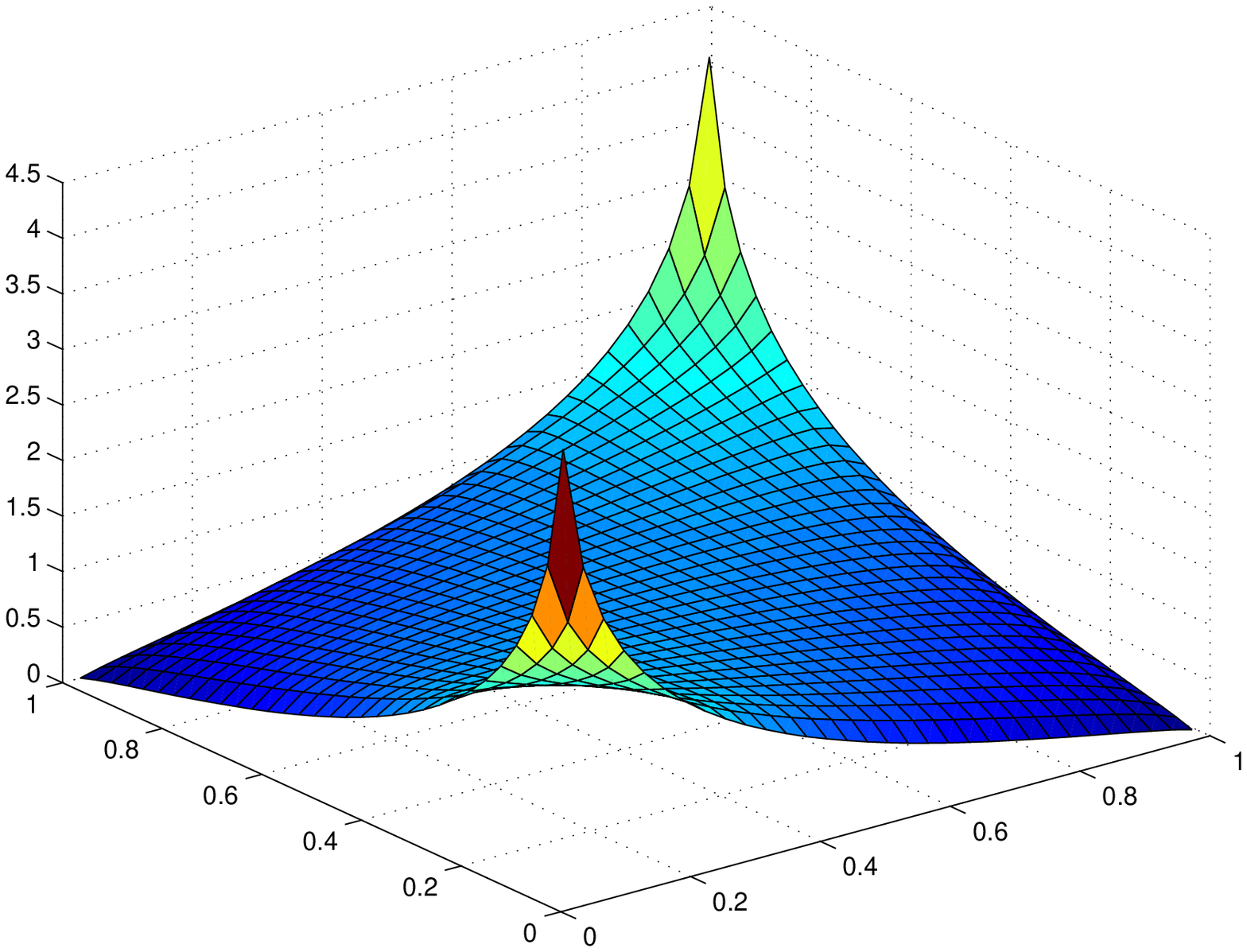}\hspace{1cm}
\includegraphics[width=4cm,height=4cm]{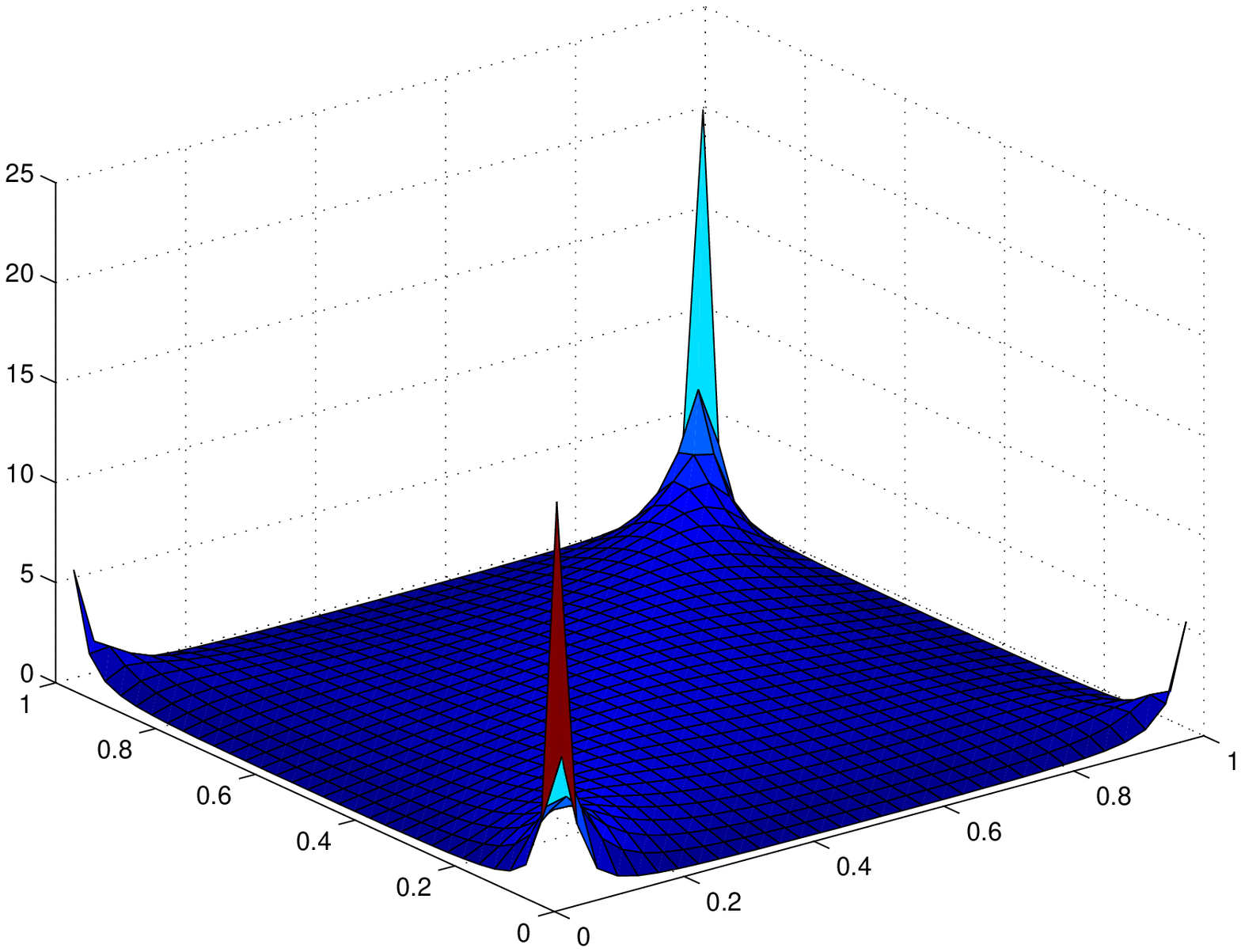}
\end{center}
\hspace{0.5cm} \caption{Kendall's tau$=0.25$. Left: Bi-dimensional
Gaussian copula density with parameter $\rho=0.4$. Right:
Bi-dimensional Student copula density with parameter
$(\rho,\nu)=(0.4,1)$. }\label{fig1}
\end{figure}

\begin{figure}[h!]
\begin{center}
\includegraphics[width=4cm,height=4cm]{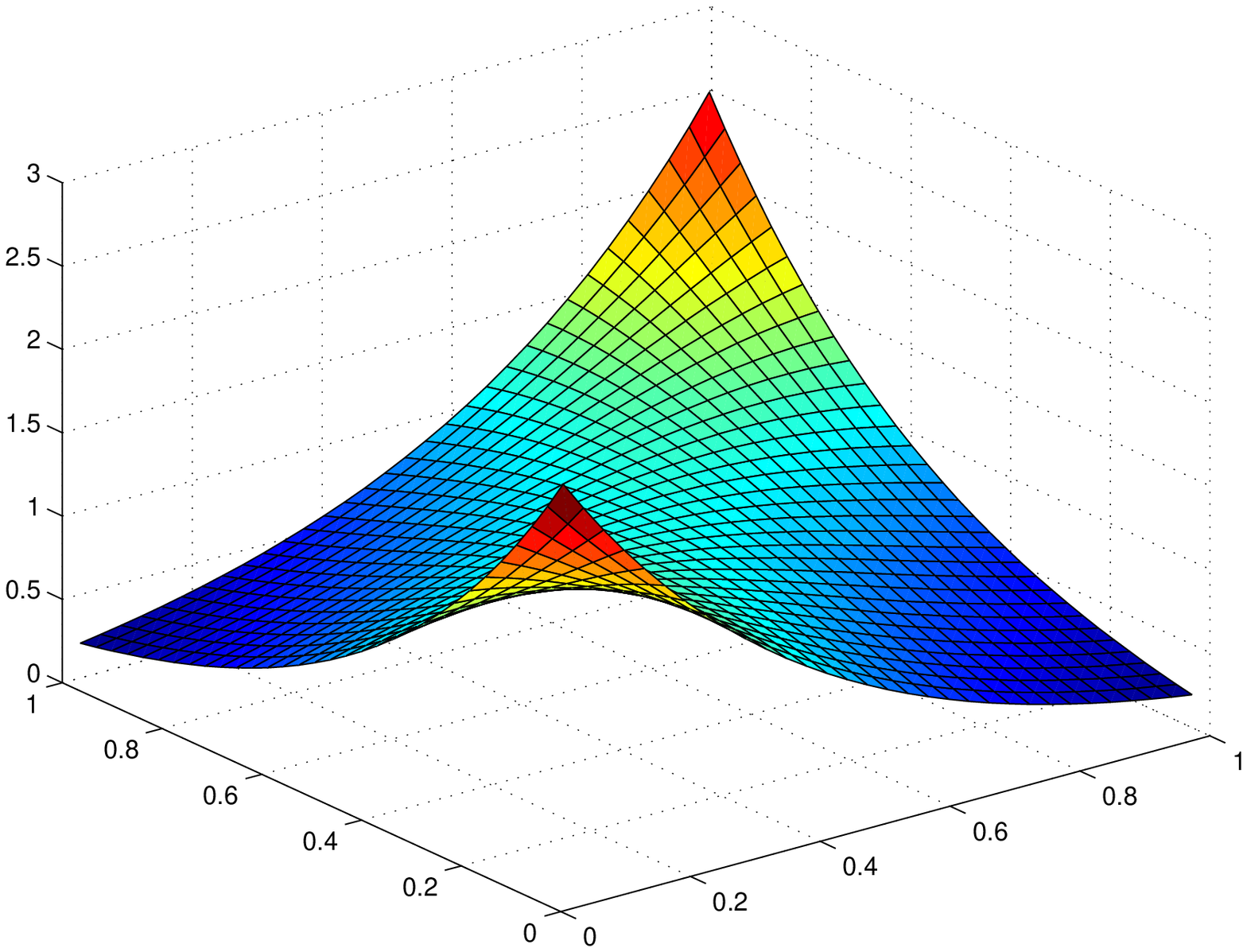}
\includegraphics[width=4cm,height=4cm]{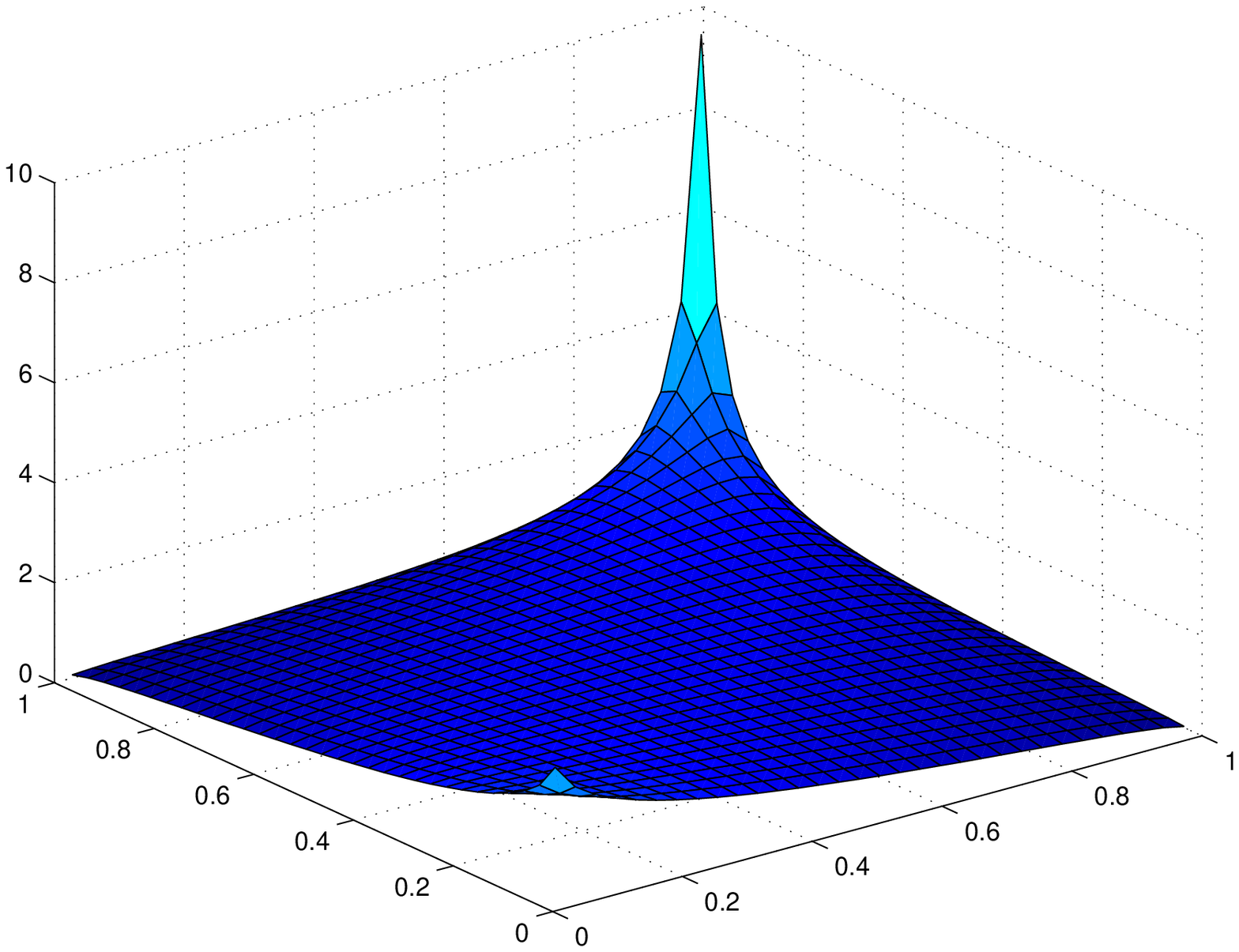}
\includegraphics[width=4cm,height=4cm]{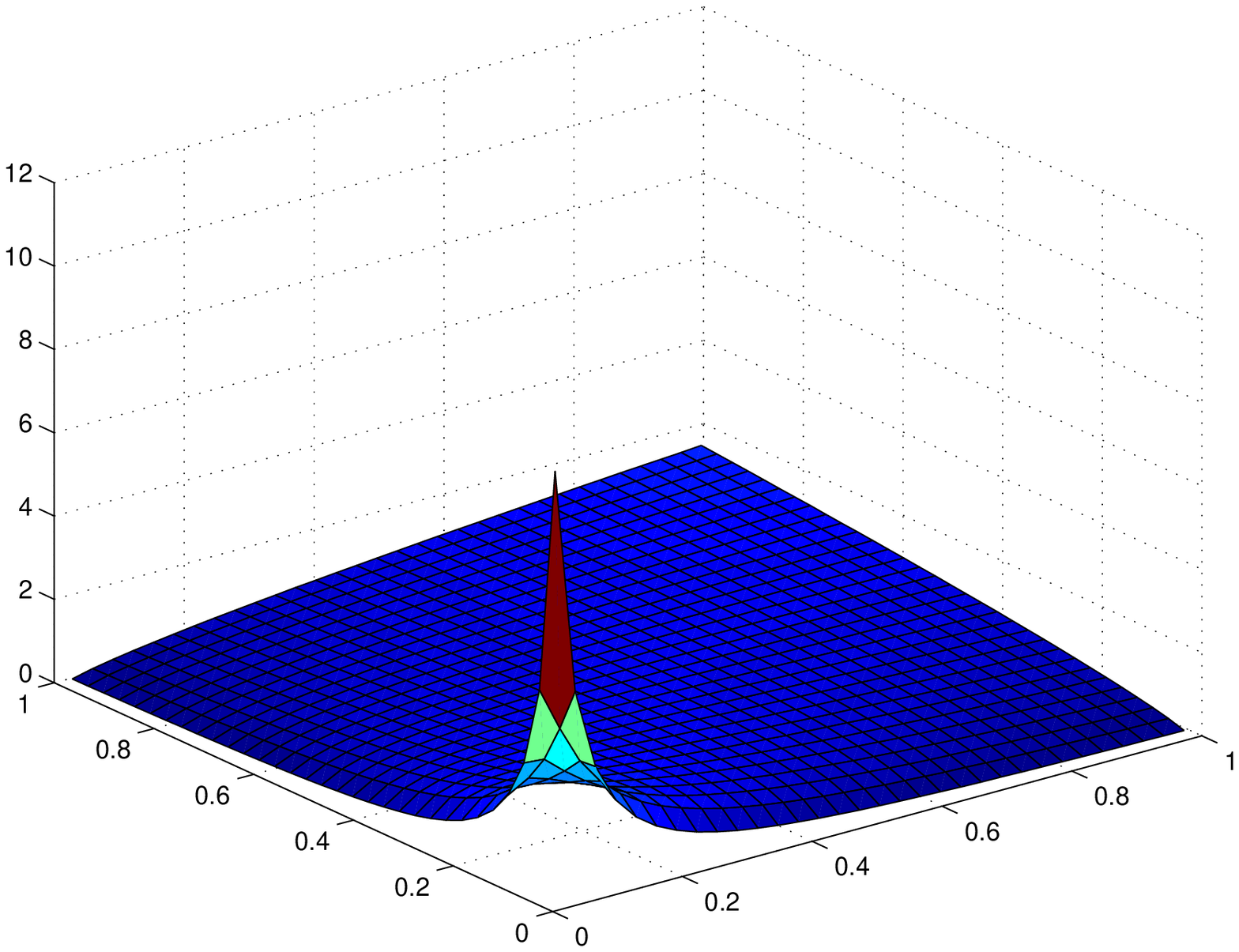}
\end{center}
\hspace{0.5cm} \caption{Kendall's tau$=0.25$. Left: Bi-dimensional
Frank copula density with parameter $\theta=2.5$. Center:
Bi-dimensional Gumbel copula density with parameter $\theta=1.33$.
Right: Bidimensional Clayton copula density with parameter
$\theta=0.66$. }\label{fig2}
\end{figure}

Since many parametric copula models are now available, the crucial
choice for the practitioner is to identify the model which is
well-adapted to  data at hand. Many goodness-of-fit tests are
proposed in the literature. \citep{Genest/Remillard/Beaudoin:2008}
give an excellent review and propose a detailed empirical study for
different tests: we refer to this paper for any supplementary
references. Roughly speaking, they study procedures based on
empirical processes. Among others, they deal with rank-based
versions of the Cram\'er-von-Mises and Kolmogorov-Smirnov
statistics. They also consider test based on Kendall's transform.
Basically, they restrict themselves to test statistics built from
empirical distributions (empirical copula or transform of this
latter). On a theoretical point of view, the asymptotic law under
the null of the test statistic is stated in a number of papers (see
by instance \citep{Deheuvels:1979}, \citep{Deheuvels:1981a} and
\citep{Deheuvels:1981b}). It allows in particular to derive the
critical value but generally the alternative is unspecified and the
properties on the power are empirically given from simulations.

In our paper, it is supposed that the copula $C$ admits a density
copula $c$ with respect to the Lebesgue measure. To our knowledge,
\citep{Fermanian:2005} was the first author to propose a
goodness-of-fit test based on nonparametric kernel estimations of
the density copula. In the same spirit as the papers cited above, he
derived the asymptotic law of the test statistic under the null. His
results are valid for bandwidths greater than $n^{-2/(8+d)}$ which
correspond to enough smooth copula densities.

Here, we focus on the minimax theory framework: we define the test
problem as initiated by \citep{Ingster:1982}. One of the advantages
of this point of view is to precisely define the alternative: it is
then possible to quantify the risk associated with the test problem
as the sum of the first type error and the second type of error.
Since this risk measure provides a quality criterion, it is then
possible to compare the test procedures. Indeed, the alternative
$H_1(v_n)$ is defined from  a positive quantity $v_n$ measuring the
distance between the null and the latter. Obviously, the larger is
this separating distance, the easier is the decision. The aim of the
minimax theory is to determine the larger alternative for which the
decision remains feasible. Solving {\bf the lower bound problem} is
equivalent to exhibit the faster separating rate $v_n$ such that the
risk is bounded from below by a given positive constant $\alpha$:
this rate is called {\bf the minimax rate of testing}. Next, {\bf
the upper bound problem} has to be solved exhibiting a test
procedure whose  risk is bounded from above by a given $\alpha$,
that is, the statistic test allows to distinguish the null from
$H_1(v_n)$, where $v_n$ is the minimax rate.

In the white noise model or in the density model,  the
goodness-of-fit problem (stands as explained above) was solved for
different regularity classes (H\"older or Sobolev or Besov)
associated with various geometries: pointwise, quadratic and
supremum norm. For fixed smoothness of the unknown density ({\bf
minimax context}), there is a rich literature summed-up in
\citep{Ingster:1993} and in \citep{Ingster/Suslina:2002}. Optimal
test procedures include orthogonal projections, kernel estimates or
$\chi^2$ procedures. Goodness-of-fit tests with alternatives of
variable smoothness into some given interval ({\bf adaptive
context}) were introduced by \citep{Spokoiny:1996} for the
$L_2$-norm
 in the Gaussian white noise model and generalized by
\citep{Spokoiny:1998} to  $L_p$-norms. \citep{Ingster:2000} proved
that a collection of $\chi^2$ tests attains the adaptive rates of
goodness-of-fit tests in $L_2$-norm as well as for the density
model.

For sake of simplicity, we restrict ourselves to bi-dimensional data
but there is no theoretical obstacle to generalize our results to
higher dimensions. Suppose that we observe $n$ i.i.d. copies
$(X_i,Y_i)_{i\in {\cal I}}$ where ${\cal I}=\{ 1,\ldots,n\}$ of
$(X,Y)$. The random vector $(X,Y)$ is drawn from the distribution
function $H$ expressed through the copula $C$. Moreover, it is
assumed that $C$ has  a copula density $c$ with respect to the
Lebesgue measure  and  $F$ and $G$ stand for the cdf's of $X$ and
$Y$ respectively. From $(X_i,Y_i)_{i\in {\cal I}}$, we are
interested in studying the goodness-of-fit problem when the null is
a composite hypothesis $H_0:c\in {\cal C}_{\Lambda}$ for a general
class ${\cal C}_{\Lambda}$ of parametrical copula densities. Since
the alternative is defined from  the quadratic distance, we propose
a goodness-of-fit test based on wavelet estimation of an integrated
functional of the copula density. Indeed,
\citep{Genest/Masiello/Tribouley:2008} and
\citep{Autin/LePennec/Tribouley:2008} show that the wavelet methods
are an efficient tool to estimate the copula densities since these
latter have very specifics behaviors. Unfortunately no direct
observations $(F(X_i),G(Y_i))$ for $i \in {\cal I}$  are available
since $F$ and $G$ are unknown, the test statistic is then built with
pseudo-observations $(\widehat{F}(X_i),\widehat{G}(Y_i))_{i\in {\cal
I}}$: as usual in the copula context, the quantities of interest are
rank-based statistics. We provide an auto-driven test procedure and
we produce its rate when the alternative contains a regular
constraint: since the procedure is based on wavelet methods, the
linked functional classes are the Besov classes $B_{s,p,q}$. We give
results for $p\geq 2$ ({\bf dense case}) and $s\geq 1/2$. The
constraint $s\geq 1/2$ is due to the fact that pseudo-data are used
and then a minimal regularity is required in order to pay no
attention to substitute the direct data with the ranked data.
Observe that \citep{Kerkyacharian/Picard:2004} have the same
constraint in the univariate regression model when the design is
random with unknown distribution. Next, we prove that our procedure
is minimax (and adaptive) optimal by exhibiting the minimax adaptive
rate. This one looks like the minimax rate but an extra $\log\log$
term appears: we prove that this loss is the price to paid for
adaptivity. To our knowledge, the proof of the adaptive lower bound
in the multivariate  density model when the null is composite has
never been clearly written.

\begin{figure}[h!]
\begin{center}
\includegraphics[width=4cm,height=4cm]{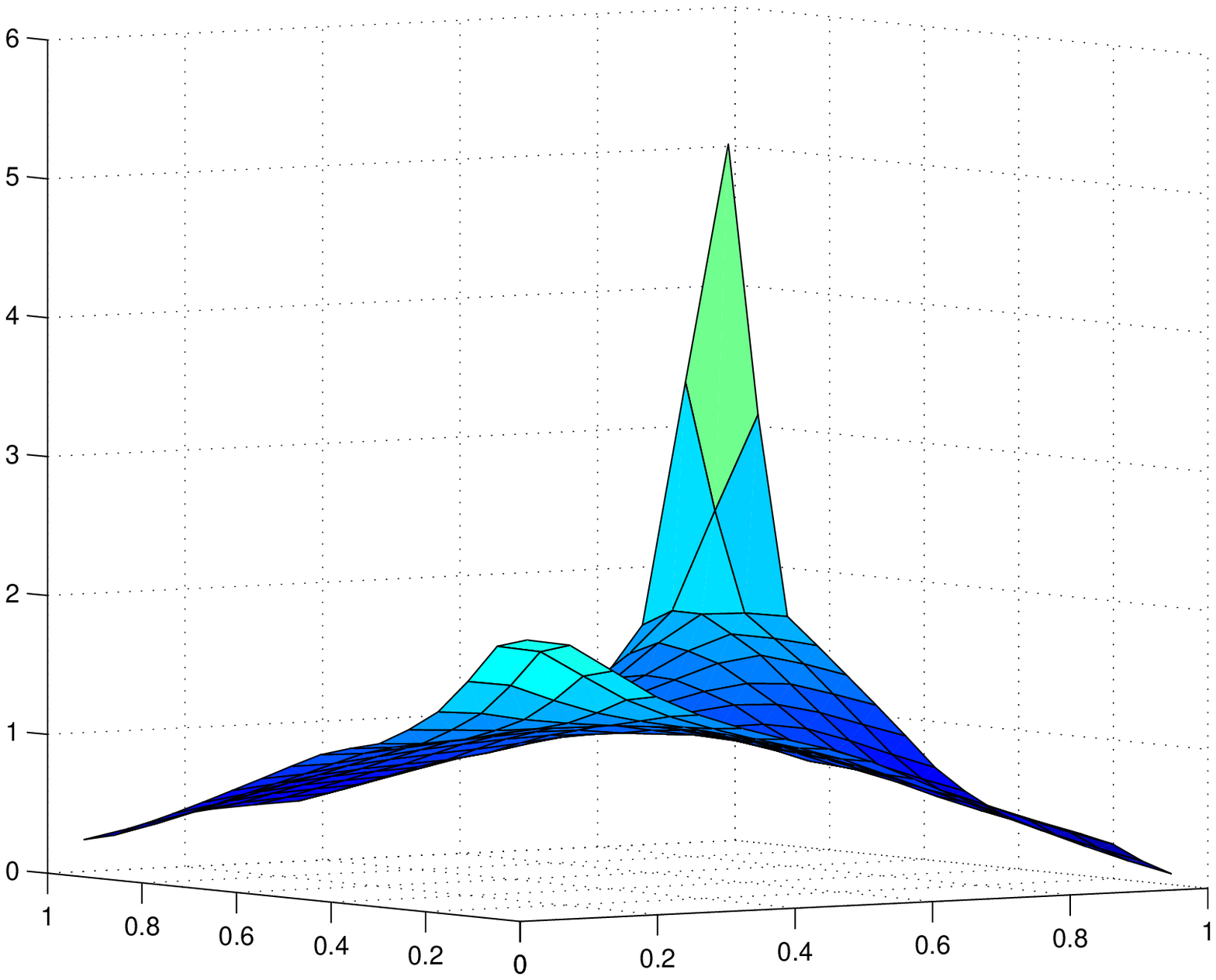}
\includegraphics[width=4cm,height=4cm]{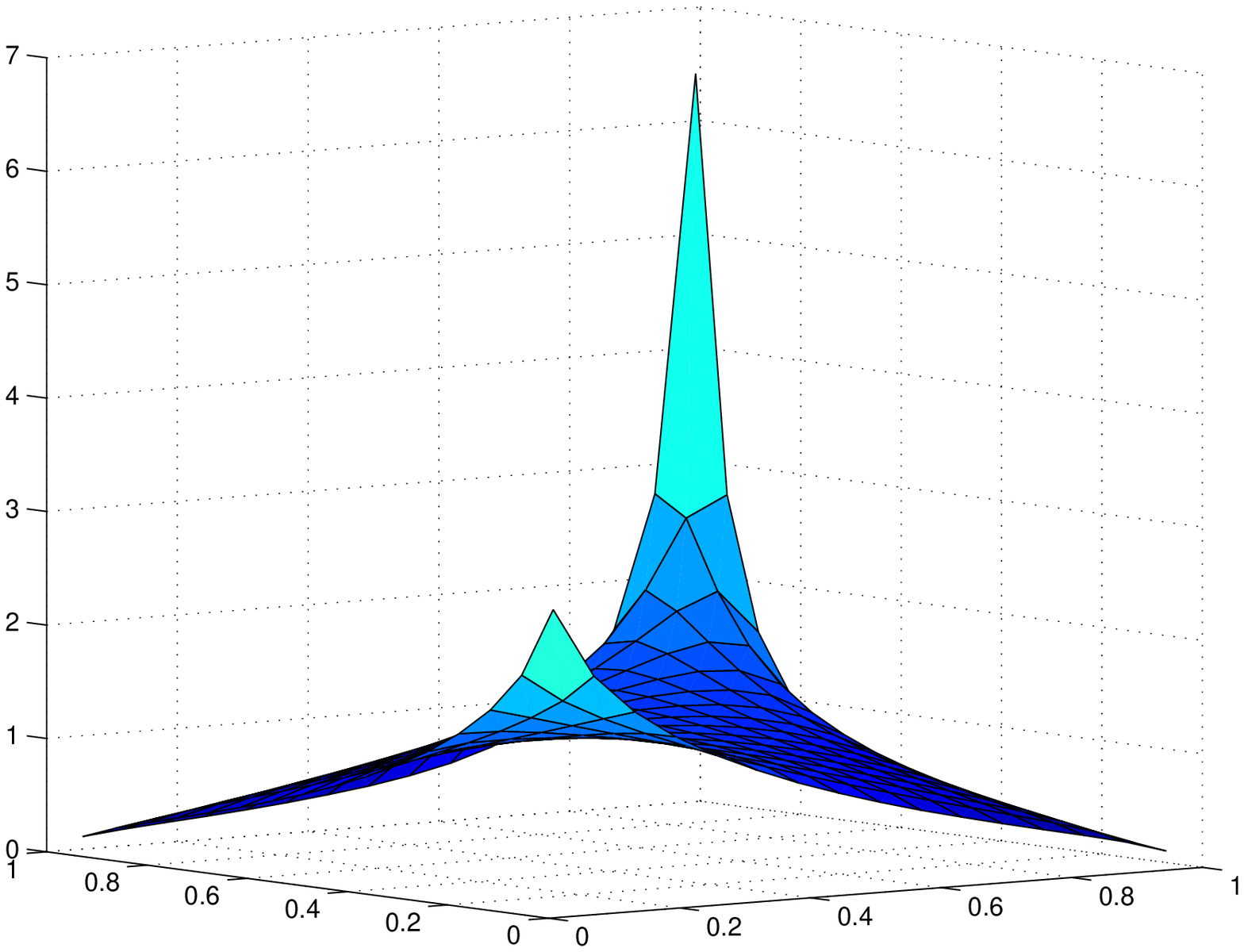}
\includegraphics[width=4cm,height=4cm]{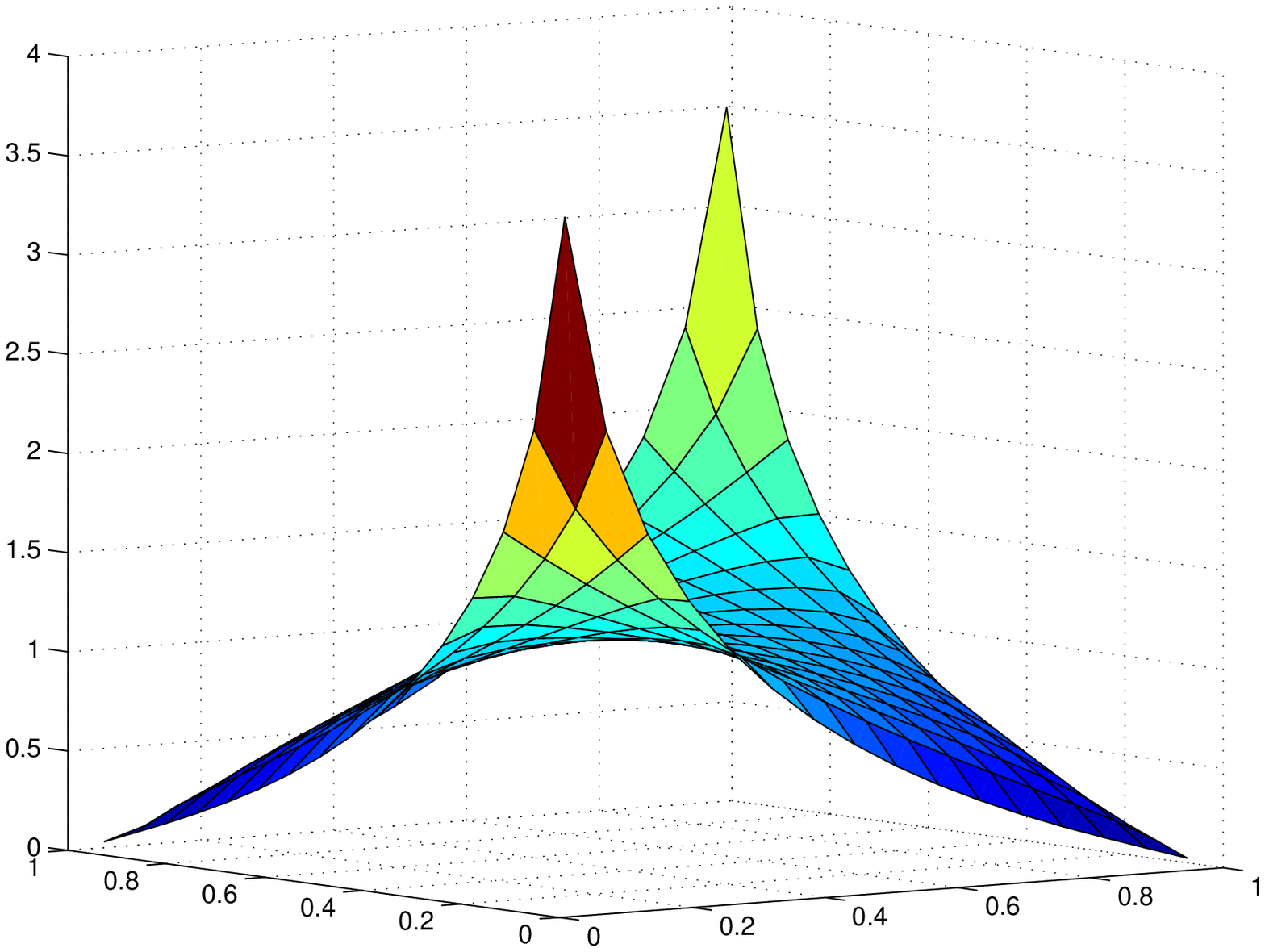}
\end{center}
\hspace{0.5cm} \caption{Left: Thresholded wavelet estimator for the
copula density of $(LOSS,ALEA)$ as given in
\citep{Autin/LePennec/Tribouley:2008}. Center: Gumbel copula density
with parameter $\theta=1.45$. Right: Gaussien copula density with
parameter $\rho=0.48$. }\label{fig3}
\end{figure}

Next, we allocate a part to empirical studies. Simulation allows us
to show that, when the theoretical framework is respected, the power
qualities of our test procedures are good. We choose to make
simulations starting from the parametrical copula families presented
at the beginning of the introduction and which are the more common
for applications. We compare our simulation results with those of
\citep{Genest/Remillard/Beaudoin:2008}. Then, we study a very well
known sample of real life data of \citep{Frees/Valdez:1998}
consisting of the indemnity payment (LOSS) and the allocated loss
adjustment expense (ALAE) for 1500 general liability claims. The
most popular model for the copula is a Gumbel copula model with
parameter $\theta=1.45$ (which may be estimated by inverting the
Kendall's tau) given in Figure \ref{fig3}. Among other results, it
is empirically shown that the Gumbel and the Gaussian copula models
are acceptable
 while Student, Clayton or Frank models are rejected. Figure \ref{fig3} gives a wavelet estimator of the
copula density of $(LOSS,ALEA)$ by the method explained in
\citep{Autin/LePennec/Tribouley:2008}. Visually, fitting the unknown
copula with the Gumbel model seems indeed to be the most
appropriated.

The paper is organized as follows. In Section 2, we first provide a
general description of orthonormal wavelet bases, focusing on the
mathematical properties that are essential to the construction of
the statistics that we consider. In Section 3, we provide the
inference procedures: first, we explain how to estimate the square
$L_2$-norm of the copula density
 and next we derive the procedure of goodness-of-fit. The
theoretical part is exposed in Section 4: first, we state very
precisely the test problem under consideration; we define the
criterion allowing to measure the quality of test procedures and
define the separating minimax rate. In Section 5, the main results
are stated: our test procedure is shown to be optimal in the sense
defined in the previous section. Section 6 is devoted to practical
results with both simulated and real data. We conclude these parts with a
discussion in Section 7. The proof of the upper bound is given in
Section 8 while the proof of the lower bound is given in Section 9.
Finally, all  technical or computational
lemmas which are not essential to understand the main proofs, are postponed in appendices.


\section{Wavelet Setting}
\subsection{Wavelet expansion}
In the univariate case, we consider a wavelet basis of $L_2([0,1])$
(see \citep{Cohen/Daubechies/Vial:1993}). Let $\phi$ be the scaling
function and let $\psi$ be the same notation for the associated
wavelet function and its usual modifications near the frontiers $0$
and $1$. They are chosen compactly supported on $[0,L]$, $L>0$. Let
$j$ in $\BBn$, $k_1$ in $\BBz$ and for any univariate function
$\Phi$, set $\Phi_{j,k_1}(\cdot)=2^{j/2}\Phi(2^j\cdot-k_1)$. In the
sequel, we use wavelet expansions for bivariate functions and we
keep the same notation as for the univariate case. Then, a bivariate
wavelet basis is built as follows:
\begin{eqnarray*}
\phi_{j,k}(x,y)=\phi_{j,k_1}(x)\phi_{j,k_2}(y),&\quad&
\psi_{j,k}^{(1)}(x,y)=\phi_{j,k_1}(x)\psi_{j,k_2}(y), \\
\psi_{j,k}^{(2)}(x,y)=\psi_{j,k_1}(x)\phi_{j,k_2}(y),&\quad&
\psi_{j,k}^{(3)}(x,y)=\psi_{j,k_1}(x)\psi_{j,k_2}(y),
\end{eqnarray*}
where the subscript $k=(k_1,k_2)$ indicates the number of components
of the functions $\phi_{j,k}$ and $\psi_{j,k}$. For a given $j\in
\BBn$, the set $$ \{\phi _{j,k},\psi _{\ell,k^\prime}^\epsilon, \;
\ell\geq j,(k,k^\prime)\in \BBz^{2}\times \BBz^2,\epsilon =1,2,3\}$$
is an orthonormal basis of $L^2([0,1]^2)$ and the expansion of any
real bivariate function $\Phi$ in $L^2([0,1]^2)$ is given by:
\[
\Phi(x,y)=\sum_{k\in \BBz^2}A _{j,k}\phi _{j,k}(x,y)+\sum_{\ell=
j}^\infty\sum_{k\in \BBz^2}\sum_{ \epsilon =1,2,3}B
_{\ell,k}^\epsilon \psi _{\ell,k}^\epsilon(x,y),
\]
where the scaling coefficients and the wavelet coefficients are
\[
\forall j\in \BBn,\forall k\in \BBz^2,\quad
A_{j,k}=\int_{[0,1]^2}\Phi \phi_{j,k},\qquad
B_{j,k}^\epsilon=\int_{[0,1]^2}\Phi \psi_{j,k}^\epsilon.
\]
The Parseval Equality immediately leads to the expansion of the
square $L_2$-norm of the function $\Phi$:
\begin{eqnarray}\label{decomposition}
\int \Phi^2&=&T_{j}+B_{j},
\end{eqnarray}
where the trend and the detail terms are respectively:
\begin{eqnarray}\label{TB}
T_{j}=\sum_{k\in \BBz^2}\left( A_{j,k}\right)^2 \mbox{ {\rm and} }
B_{j}=\sum_{\ell=j}^\infty\;\sum_{k\in \BBz^2}\;\sum_{\epsilon=1}^3(
B_{\ell,k}^\epsilon)^2.\end{eqnarray} Notice that, since the support
of $\Phi$ is  $[0,1]^2$, the sum over the indices $k$ is finite:
there are no more than $(2^j+L)^2$ terms in the sum (recall that $L$
is the length of the support of $\phi$). In order to simplify the
notations, the bounds of variation of $k$ and $\epsilon$ in
expansion of any $\Phi$, are omitted in the sequel.

\subsection{Besov Bodies and Besov spaces}
Dealing with wavelet expansions, it is natural to
consider Besov bodies as functional spaces since they are
characterized in term of wavelet coefficients as follows.
\begin{definition}\label{BesovBody}
For any $s>0, \; p\geq 1$ and any radius $M>0$,
a $d-$varied function $\Phi$ belongs to the ball $b_{s,p,\infty}(M)$ of the Besov
body $b_{s,p,\infty}$
if and only if
its sequence of wavelet coefficients $B_{j,k}^\epsilon$ satisfies
\begin{eqnarray*}
\forall j\in \BBn,\;\sum_{k\in
\BBz^2}\sum_{\epsilon=1}^3|B_{j,k}^\epsilon|^p <M2^{-j(s+d/2-d/p)p}.
\end{eqnarray*}
\end{definition}
The Besov body $b_{s,p,\infty}$ coincides with the more standard
Besov space $\mathcal{B}_{s,p,\infty}$ when there exists an integer
$N$ strictly larger than $s$ and such that the $q-$th moment of the
wavelet $\psi$ vanishes for any $q=0,\ldots, N-1$. It is possible to
build univariate wavelets whose support is included in $[0, 2N -1]$
satisfying  this property for any choice of $N$ (see the
Daubechies wavelets).

In the sequel, we
need to bound the detail term $B_j$ defined in (\ref{TB}). We use
the  following inequality
\begin{eqnarray*}
\forall j\in \BBn,\quad B_j\leq
\sum_{\ell=j}^\infty\;\left(\sum_{k\in
\Z^2}\;\sum_{\epsilon=1}^3\left|B_{\ell,k}^\epsilon\right|^p\right)^{2/p}\left(K\,2^{2j}\right)^{1-2/p},
\end{eqnarray*}
where $K$ is a positive constant depending on the supports of $\Phi$
and $\psi$.
Assuming that the function $\Phi$ belongs to
$b_{s,p,\infty}(M)$ with $s,p$ and $M$ as in Definition \ref{BesovBody}, the following
inequality
holds
\begin{eqnarray}\label{reg}
\forall j\in \BBn, \quad B_j&\leq& {\tilde K}\;2^{-2js},
\end{eqnarray}
where ${\tilde K}$  is a positive constant depending on the supports
of $\Phi$, $\psi$ and on the radius $M$. When $\Phi$ is a copula
density, ${\tilde K}=M^{2/p}\left( 3(L+1)^2 \right)^{1-2/p}$.

\section{Statistical Procedures}\label{procedure}
Assuming that the copula density $c$
belongs to $L_2([0,1]^2)$, we first explain the procedure to
estimate the square $L_2-$norm of $c$
$$
\theta=\,\|c\|^2:=\,\int_{[0,1]^2} c^2,
$$
which is used to define the alternative of the goodness-of-fit test.
The statistical methods depend on parameters (the level $j$ for the
estimation procedure and $j$ and the critical value $t_j$ for the
test procedure) which are discussed and determined in an optimal way
in Section \ref{mainR}.

It is fundamental to Notice that, for any bivariate function $\Phi$,
one has
\begin{eqnarray}\label{passage_c_h}
\E_c\left[\Phi(U,V)\right]&= &\E_h\left[\Phi(F(X),G(Y))\right],
\end{eqnarray}
where $h$ stands for the joint density  of $(X,Y)$. This means in
particular that the wavelet coefficients $
\{c_{j,k},c_{\ell,k}^\epsilon, \; \ell\geq j,k\in
\BBz^2,\epsilon=1,2,3\} $ of the copula density $c$ on the wavelet
basis $$ \{\phi _{j,k},\psi_{\ell,k}^\epsilon, \; \ell\geq j,k\in
\BBz^2,\epsilon=1,2,3\}$$ are equal to the coefficients of the joint
density $h$ on the warped wavelet family
$$ \{\phi_{j,k}(F(\cdot), G(\cdot)),\psi
_{\ell,k}^\epsilon(F(\cdot), G(\cdot)), \; \ell\geq j,k\in
\BBz^2,\epsilon =1,2,3\}.$$ The statistical procedures are based on
the wavelet expansion of the copula density $c$, for which the
wavelet coefficients have to be estimated.
\subsection{Procedures to estimate $\theta$} \label{ESTprocedure}
Let $J$ be a subset of $\BBn$ and consider a given $j$ in $J$.
Motivated by the wavelet expansion (\ref{decomposition}), we propose
to estimate $\theta$ with an estimator of the trend $T_{j}$ omitting
the detail term $B_j$. Using the orthonormality property of the wavelet basis,
it leads to estimate the square of the
coefficients of the copula density on the scaling function. As
usual, a $U-$statistic associated with the empirical coefficients is
used in order to remove the bias terms. Due to (\ref{passage_c_h}),
we first consider the following family of statistics $\{\widehat{
T_j},j\in J\}$ defined by
\begin{eqnarray*} \label{That}
\widehat{ T_j}&=&\sum_{k} \widehat{ \theta_{j,k}},
\end{eqnarray*}
where $\widehat{ \theta_{j,k}}$ is the following $U-$statistic
\begin{eqnarray*} \widehat{ \theta_{j,k}} &=&
\frac{1}{n(n-1)} \mathop{\sum^n_{i_1, i_2 =1}}_{ i_1 \neq i_2
}\phi_{j,k}\left(F(X_{i_1}), G(Y_{i_1})\right) \phi_{j,k}\left(
F(X_{i_2}), G(Y_{i_2})\right).
\end{eqnarray*}
Since no direct observation $(F(X_i),G(Y_i))$
is usually
available, it is  replaced  in $\widehat{ \theta_{j,k}}$ by the
pseudo observation $(\widehat {F}(X_i), \widehat {G}(Y_i))$, where
$\widehat {F}, \widehat {G}$ denote some estimator of the margins.
To preserve the independence given by the observations, we split the
initial sample $(X_i,Y_i)_{i \in {\cal I}}$ into disjoint samples
$(X_i,Y_i)_{i \in {\cal I}_1}$ and $(X_i,Y_i)_{i \in {\cal I}_2}$
with ${\cal I}_2 \cup {\cal I}_1={\cal I}, \; {\cal I}_2 \cap {\cal
I}_1=\emptyset$, and whose size is $n_1$ and $n_2$ respectively. The
sub-sample with indices in ${\cal I}_1$ is used to estimate the
marginal distributions and the second one with indices in ${\cal
I}_2$ is devoted to the computation of the $U$-statistic. We
consider the usual empirical distribution functions:
\begin{eqnarray*}
\widehat{F}(x)=\frac{1}{n_1} \mathop{\sum_{i \in {\cal I}_1}} \I_{\{X_i \leq x\}}&\mbox{
and }&
\widehat{G}(y)=\frac{1}{n_1}\mathop{\sum_{i \in {\cal I}_1}}
\I_{\{Y_i \leq y\}} .
\end{eqnarray*} It leads
to the family  $\{\widetilde{T_j},j\in J\}$ of estimators of
$\theta$
\begin{eqnarray*}\label{Ttilde}
\widetilde{T_{j}}&=&\sum_{k}\widetilde{\theta_{j,k}},
\end{eqnarray*} with \begin{eqnarray*}
\widetilde{\theta_{j,k}}&=&\frac{1}{n_2(n_2-1)}\mathop{\sum_{i_1,
i_2\in {\cal I}_2}}_{ i_1 \neq i_2 } \phi_{j,k}\left(\frac{R_{i_1}}{n_1}
,\frac{S_{i_1}}{n_1} \right) \phi_{j,k}\left(\frac{R_{i_2}}{n_1},
\frac{S_{i_2}}{n_1} \right),\nonumber
\end{eqnarray*}
where $R_p=n_1 \widehat{F}(X_p)$ and $S_p = n_1 \widehat{G}(Y_p)$, $p \in {\cal I}_1$,
could
be viewed as estimates of the rank statistics of  $X_p$ and $Y_p$
respectively.

\subsection{Test Procedures} \label{TESTprocedure}
In this part, we consider a family of known bivariate copula
densities ${\cal
C}_\Lambda=\left\{c_\lambda,\lambda\in\Lambda\right\}$
indexed
by a parameter $\lambda$ varying in a given set $\Lambda \subset
\BBr^{d_{\Lambda}}$, $d_{\Lambda} \in \BBn^*$. From the observations
$(X_i,Y_i)_{i \in {\cal I}}$, our aim is to test the goodness-of-fit
between any $c_\lambda$ and a copula density $c$, which is {\it
enough distant} in the $L_2$-norm, from the parametric family ${\cal
C}_\Lambda$. Acting as in paragraph \ref{ESTprocedure},
we estimate the square $L_2$-norm between $c$ and a fixed element $c_\lambda$ lying in
the family
${\cal C}_\Lambda$ by
\begin{eqnarray}\label{Ttilde}
\widetilde {T_j}(\lambda)&=&\sum_{k}\widetilde
{\theta_{j,k}}(\lambda),
\end{eqnarray} for \begin{eqnarray*}
\widetilde
{\theta_{j,k}}(\lambda)&=&\frac{1}{n_2(n_2-1)}\mathop{\sum_{i_1,
i_2\in {\cal I}_2}}_{ i_1 \neq i_2 } \left(\phi_{j,k}\left(\frac{R_{i_1}}{n_1}
,\frac{S_{i_1}}{n_1} \right)- c_{j,k}(\lambda)
\right)\\&&\hspace{3cm} \times
\left(\phi_{j,k}\left(\frac{R_{i_2}}{n_1}, \frac{S_{i_2}}{n_1}
\right)- c_{j,k}(\lambda) \right)\nonumber,
\end{eqnarray*}
where  $ \{c_{j,k}(\lambda),k\in\BBz^2,j\in \BBn\}$ denote the known
scaling coefficients of the target copula density $c_\lambda$.
Notice that, if  direct observations $(F(X_i),G(Y_i))_{i \in
{\cal I}}$ would be available, the appropriate test statistic $\widehat{
T}_j (\lambda)$ would be
\begin{eqnarray*} \label{ThatLambda}
\widehat{ T}_j (\lambda)& =& \sum_k \widehat{\theta}_{j,k}(\lambda),
\end{eqnarray*}
where
\begin{eqnarray*}
\widehat
{\theta_{j,k}}(\lambda)&=&\frac{1}{n_2(n_2-1)}\mathop{\sum_{i_1,
i_2\in {\cal I}_2}}_{ i_1 \neq i_2 }
\left(\phi_{j,k}\left(F(X_{i_1}),G(Y_{i_1}) \right)-
c_{j,k}(\lambda) \right)\\&&\hspace{3cm} \times
\left(\phi_{j,k}\left(F(X_{i_2}),G(Y_{i_2}) \right)-
c_{j,k}(\lambda) \right)\nonumber.
\end{eqnarray*}
Now we are ready to build the test procedures. Let
us give
a set of indices $J$ and a set of critical values $\{t_j, \; j \in
J\}$ and define $\{D_j^\Lambda,j\in J\}$, the family
of test statistics
\begin{eqnarray*}\label{testopti}
D_j^\Lambda&=&  \I_{ \displaystyle{\inf_{\lambda \in \Lambda}}
\widetilde {T_j}(\lambda) > t_j} \quad ,\end{eqnarray*} allowing to
test if $c$ belongs to the parametric family ${\cal
C}_\Lambda=\{c_\lambda,\lambda\in\Lambda\}$. Notice that
$\Lambda=\{\lambda_0\}$ leads to the single null hypothesis
$H_0:c=c_{\lambda_0}$. We are also interested in building
auto-driven procedures by considering all the tests in the family
\begin{eqnarray}
D_\Lambda&=& \max_{j\in J}D_j^\Lambda=  \I_{ \displaystyle{ \max_{j\in J}
(\inf_{\lambda \in \Lambda} \widetilde {T_j}(\lambda)- t_j  )}>
0}\;. \label{testAdap}
\end{eqnarray}
The sequence of parameters $t_j$ of the method are determined in an
optimal way in Section 5. We explain in Section 4 what ``optimal way'' means
in giving a presentation of the minimax theory for our
framework.
\section{Minimax Theory}
We adopt the  minimax point of view to solve the problem of
hypothesis testing, initiated by \citep{Ingster:1982} in Gaussian
white noise. A review of results obtained in  problems of minimax
hypothesis testing is available in \citep{Ingster:1993} and
\citep{Ingster/Suslina:2002}. Let us describe this approach.
\subsection{Minimax hypothesis testing Problem}
As in the previous section, we consider ${\mathcal
C}_{\Lambda}=\{ c_\lambda, \lambda \in \Lambda\}$  a  given functional class of copula densities. For any given
$\tau=(s,p,M)$, with  $s>0,p\geq 1,M>0$, the following statistical problem of hypothesis testing is considered,

\begin{eqnarray}
H_0 : c=c_\lambda \in {\cal C}_{\Lambda}\quad  &\mbox{ {\rm against}
}&\quad  H_1 :  c \in \Gamma(v_n(\tau)),\label{pbtest}\end{eqnarray}
with
\begin{eqnarray*}\Gamma(v_n(\tau))= b_{s,p,\infty}(M)\cap \left\{c : \inf_{c_\lambda
\in {\cal C}_{\Lambda}}\|c -c_\lambda\| \geq v_n(\tau) \right\},
\end{eqnarray*}
 where $b_{s,p,\infty}(M)$ is the ball of
radius $M$ of the Besov body $b_{s,p,\infty}$ defined in Definition
\ref{BesovBody} and $v_n(\tau)$ is a sequence of positive numbers,
depending on $\tau$ and  decreasing to zero as $n$ goes to infinity.
Recall that $\|g\|$ denotes the $L_2$-norm of any function $g $ in
$L_2([0,1]^2)$. Observe that the functional class
$\Gamma(v_n(\tau))$, which determines the alternative $H_1$, is
characterized by three parameters: the regularity class
$b_{s,p,\infty}$ where the copula density is supposed to belong, the
$L_2$-norm which is the geometrical tool measuring the distance
between both hypotheses, and the sequence $v_n(\tau)$.

According to the principle of the minimaxity,  the regularity space
and the loss function are chosen by the statistician. Notice that
the parameter  $\tau$ could be known or unknown. Obviously, our aim
is to consider tests which are able to detect  alternatives defined
with sequences $v_n(\tau)$ as small as possible. It can be shown
(\citep{Ingster:1993}) that $v_n(\tau)$ cannot be chosen in an
arbitrary way: indeed, if $v_n(\tau)$ is too small, then $H_0$ and
$H_1$ cannot be distinguished with a given error $\alpha \in (0,1)$.
Therefore, solving hypothesis testing problems via the minimax
approach consists in determining the smallest sequence $v_n(\tau)$
for which such a test is still possible and to indicate the
corresponding test functions. The smallest sequence $v_n(\tau)$ is
called the minimax rate of testing. Let  $D_n$ be
 a {\it test statistic} i.e. an arbitrary function with possible values
$0,1$, measurable with respect to $(X_i,Y_i)_{i \in {\cal I}}$ and such
that we accept $H_0$ if $D _n=0 $ and we reject it if $D_n=1$.

\vspace{0.2cm}

\begin{definition} Assuming $\tau$ to be known, the sequence $v_n(\tau)$ is the minimax rate
of testing $H_0$ versus $H_1$ if  relations (\ref{Bi}) and
(\ref{BS1}) are fulfilled:
\begin{itemize}
\item for any given $\alpha_1 \in (0,1)$, there exists $a>0$ such that
\begin{eqnarray} && \lim_{n\rightarrow +\infty}\inf_{{D}_n} \left(
 \sup_{c_\lambda \in {\cal C}_{\Lambda}} \P_{\lambda} ({D}_n=1)
  +\sup_{c \in \Gamma(a\;v_n(\tau))} \P_c ({D}_n=0) \right) \geq \alpha_1 , \!\!\!\! \label{Bi}
\end{eqnarray}
where the infimum is taken over any test statistic $D_n$,
\item there exists a sequence of test statistics $(D_n^\star)_n$  for which for any given $\alpha_2 $ in
$(0,1)$, it exists $A>0$ such that
\begin{eqnarray}&&
 \lim_{n\rightarrow +\infty} \left(\sup_{c_\lambda \in {\cal C}_{\Lambda}} \P_{\lambda} (D_n^\star=1)+\sup_{c \in \Gamma(A\;v_n(\tau))} \P_c
(D_n^\star=0) \right)\leq \alpha_2,
 \label{BS1}
\end{eqnarray}
\end{itemize}
where $\P_c$, respectively $\P_{\lambda}$ denotes the distribution
function associated with the copula density $c$, respectively
with $c_\lambda$.
\end{definition}
\subsection{Adaptation}
Nevertheless,  since the copula function itself is unknown, the a
priori knowledge on  $\tau$  could appear unrealistic.
Therefore, the purpose of this paper is to solve the previous
problem of test in an adaptive framework i.e. in supposing that
$\tau=(s,p,M)$ is unknown but varying in a known set ${\cal S}$.
 Comparing the adaptive case with the non-adaptive case, it has been proved in different frameworks
  that a loss of efficiency in the rate of testing is unavoidable (see for instance \citep{Spokoiny:1996}, \citep{Gayraud/Pouet:2005}). This
  loss is expressed as $t_n$, a positive constant or a sequence of positive
 numbers increasing to infinity with $n$
  (as slow as possible), which appears in the rate of
 testing $v_{nt_n^{-1}}(\tau)$. Similarly to the minimax rate of testing, we define the adaptive minimax
 rate of testing as follows.
\begin{definition}
 The sequence $v_{nt_n^{-1}}(\tau)$ is
 the adaptive minimax rate of testing if  relations (\ref{borneInf}) and (\ref{borneSup}) are satisfied
\begin{itemize}
\item
for any given $\alpha_1 \in (0,1)$, there exists $ a>0$ such that
\begin{eqnarray}&&
 \lim_{n\rightarrow +\infty}\inf_{{D}_n} \left(
 \sup_{c_\lambda \in {\cal C}_{\Lambda}} \P_{\lambda} ({D}_n=1) +
\sup_{\tau \in {\cal S}} \sup_{c \in \Gamma(a\;
v_{nt_n^{-1}}(\tau))} \P_c ( { D}_n =0) \right)\geq \alpha_1 ,
\hspace{1cm} \label{borneInf}
\end{eqnarray}
where the infimum is taken over any test statistic $D_n$,
\item
there exists a sequence of universal test statistics
 $D_n^{\star}$ (free of $\tau$) such that, for any given
 $\alpha_2$ in
 $(0,1)$, there exists $A>0$ such that
\begin{eqnarray}&&
\lim_{n\rightarrow +\infty}  \left(
 \sup_{c_\lambda \in {\cal C}_{\Lambda}} \P_{\lambda} (D_n^{\star}=1)+\sup_{\tau \in {\cal
S} }   \sup_{c \in \Gamma(A\;v_{nt_n^{-1}}(\tau))} \P_c (
D_n^{\star} =0)\right)  \leq \alpha_2 \hspace{1cm} \label{borneSup}
\end{eqnarray}
where $t_n$ is either a positive constant or a sequence of positive
numbers increasing to infinity with $n$ as slow as possible.
\end{itemize}
\end{definition}
 Notice that  relations (\ref{borneInf}) and (\ref{borneSup}) (instead
 of  relations (\ref{Bi}) and (\ref{BS1})) mean that the minimax  rate of testing
$v_n(\tau)$ is contaminated by the term $t_n$ in the adaptive
setting. Observe that the same phenomenon is observed in the
estimation problem where an extra logarithm term ${\tilde t}_n=\log(n)$ has
often (but not always) to be paid for the adaptation.

\section{Main results} \label{mainR}
In this section, we focus on test problems for which the parametric
family ${\cal C}_\Lambda$ is included in some
$b_{s_\Lambda,p_\Lambda,\infty}(M_\Lambda)$ where $s_\Lambda>0$,
$p_\Lambda\geq 1$ and $M_\Lambda >0$ are known.

Our theoretical results concern the minimax resolution of the
problem of hypothesis testing defined in (\ref{pbtest}) in an
adaptive framework.
Theorem \ref{BInf} states the result of the lower bound (see
relation (\ref{borneInf})). Then, Theorem \ref{BSup} exhibits the
rate achieved by the test procedure proposed in Section 3 (see
relation (\ref{borneSup})). Comparing the rate of our procedure with
the fastest rate given in Theorem \ref{BInf} leads to Theorem
\ref{opt} establishing the  optimality of our procedure.\\
First, let us state  the assumption which gives a control of the complexity of
${\cal C}_\Lambda$.
\begin{itemize}
\item {\bf A0:} the set $\Lambda$ is  compact  in $\BBr^{d_{\Lambda}}$ and
  $$\sup_{(x,y) \in [0,1]^2} |c_\lambda(x,y)-c_{\lambda'}(x,y) | \leq Q \|\lambda-\lambda'\|_{\BBr^{d_{\Lambda}}}^\nu, \; \forall \lambda, \; \lambda' \in \Lambda, $$ where $\nu$ is a positive real, $Q$ is a positive constant and $\| \cdot \|_{\BBr^{d_{\Lambda}}}$ denotes the Euclidean norm in $\BBr^{d_{\Lambda}}$.
\end{itemize}

\subsection{Lower Bound}\label{LB}
As it is usual for composite null hypotheses, the result of the
lower bound requires the existence of a particular density
$c_{\lambda_0} \in {\cal C}_\Lambda$ (see assumption {\bf AInf} below) in
order to construct a randomized class of functions which must be
included in the alternatives.
\begin{itemize}
\item {\bf AInf: } there exists a parameter $\lambda_0$ in  $\Lambda$ such that
$$\forall (u,v)\in[0,1]^2, \quad
c_{\lambda_0}(u,v)
>m, \; \mbox{ {\rm with }} m >0.$$
\end{itemize}

\begin{theorem} \label{BInf}
Suppose that ${\cal S}$ defined by
\begin{eqnarray}\label{S}
&&{\cal S}=\{\tau=(s,p,M) , s \geq 1/2, p \geq 2,M >0: \; s-2/p\leq
s_\Lambda-2/p_\Lambda, M_\Lambda \leq M\}\hspace{1.2 cm}
\end{eqnarray}
 is nontrivial  (see \citep{Spokoiny:1996}),  which means
that there exist $p\geq 2$, $M>0$ and $0<s_{\mbox{min}} <
s_{\mbox{max}}$ such that
$$
\forall s \in [s_{\mbox{min}}, s_{\mbox{max}}],\quad (s,p,M) \in
{\cal S}
$$
 and assume that {\bf A0} and {\bf AInf} hold. Set
$$
v_{nt_n^{-1}}(\tau)=(nt_{n}^{-1})^{-2s/(4s+2)} \mbox{ {\rm with} } \;
t_{n}=\sqrt{\log (\log(n))}.
$$ Then, it exists a positive constant
$a$ such that
\begin{eqnarray}
&&\lim_{n\rightarrow +\infty}\left( \inf_{D_n}\{ \sup_{\lambda \in
\Lambda} \P_\lambda( D_n=1) + \sup_{\tau \in {\cal S}} \sup_{c \in
\Gamma(a\; v_{nt_n^{-1}}(\tau))} \P_c ( { D}_n =0) \}\right)= 1,
\label{BI}
\end{eqnarray}
where  the infimum is taken over any test function $D_n$.
\end{theorem}

\subsection{Upper Bound} \label{UB}

 Theorem  \ref{BSup} deals with relation
(\ref{borneSup}) which holds for the test statistic $D_{\Lambda}$
defined by relation (\ref{testAdap}) as soon as the parameters of
the methods are chosen as follows.
 The set
$J=\{\lfloor j_0\rfloor,\ldots,\lfloor j_\infty \rfloor\}$ is
determined by
\begin{eqnarray}\label{j0jinfty}
&&2^{j_0} = \log (n_2)\log (n_1),\quad
2^{j_\infty}=\left(\frac{n_2}{\log(n_2)}\right)^{1/2}\wedge
\left(\frac{n_1}{\log(n_1)}\right)^{1/2-1/2q},
\end{eqnarray}
where $q$ is the order of differentiability of the scaling function
$\phi$. The critical values satisfy
\begin{eqnarray}&&
\forall j\in J,\quad t_j=3 \mu\;\frac{2^j}{n_2}\sqrt{\log\log(n_2)},
\label{seuil}
\end{eqnarray}
 where $\mu$ is a positive constant such that $\mu > \sqrt{2 K_g K_1}$, and
$K_g$ and $K_1$ are positive constants depending on $\|\phi\|_\infty$, $\|c\|_\infty$,
$\|c_\lambda\|_\infty$ and the length of the
support of $\phi$  (see Lemma \ref{LDUstat}).
%

\begin{theorem} \label{BSup} Let us choose $n_1=\pi\,\,n$ and
$n_2=(1-\pi)n$ for some $\pi$ in $(0,1)$. Assume that the scaling
function $\phi$ is continuously $q-$differentiable for
$$
q\geq
\left[1-\frac{\log\left(\frac{n_2}{\log(n_2)}\right)}
{\log\left(\frac{n_1}{\log(n_1)}\right)}\right]^{-1}.
$$
Moreover assume that any density $c$ under the alternatives or any
$c_\lambda$ under the null  are uniformly bounded. Then, the test
statistic $D_\Lambda$ defined by (\ref{testAdap}) is such that
\begin{eqnarray}
 \lim_{n_1\wedge n_2 \rightarrow ¨+ \infty} \sup_{c_\lambda \in {\cal
C}_\Lambda} \P_{\lambda} ( D_\Lambda =1) = 0 \label{PremE}.
\end{eqnarray}
Assume that {\bf A0} holds, then
there exists a positive constant $A$ such that
\begin{eqnarray}
\lim_{n_1\wedge n_2 \rightarrow +\infty} \sup_{\tau \in {\cal S}}
\sup_{c \in \Gamma (A v_{nt_n^{-1}}(\tau))}  \P_{c} ( D_\Lambda =0)
&= &0, \label{SecondE}
\end{eqnarray}
where
$$
v_{nt_n^{-1}}(\tau)=(n_2t_{n_2}^{-1})^{-2s/(4s+2)}
\mbox{ {\rm and}  } \; t_{n_2}=\sqrt{\log (\log(n_2))}.
$$
\end{theorem}
Relation (\ref{borneSup}) of the upper bound holds since both
 relations (\ref{PremE}) and (\ref{SecondE}) are satisfied. Notice also  that
relation (\ref{PremE}) indicates that the test statistic
$D_\Lambda$ is asymptotically of any level in $(0,1)$.

\subsection{Optimality}
As a corollary of Theorem \ref{BInf} and Theorem \ref{BSup}, we
obtain
\begin{theorem} \label{opt}
Under  the assumptions of Theorem \ref{BInf} and Theorem \ref{BSup},
our test procedure defined by Relation (\ref{testAdap}) is adaptive
optimal over the range of parameters $\tau \in {\cal S}$ where
${\cal S}$ is defined by  equation (\ref{S}).
\end{theorem}

\section{Practical results}
The purpose of this section is to provide several examples to
investigate the performances of the test procedure presented in
Section \ref{procedure}. This part is not exactly an illustration of
the theoretical part since it does not focus on the separating rate
between the alternative and the null hypothesis, but it is devoted
to the study of our test procedure from  a risk point of view.
Note also that we do not use exactly the theoretical procedure
described in the previous section. As usual for practical purpose,
we replace theoretical quantities by more adapted quantities
 obtained with resampling methods. In the first part, we fix the
test level $\alpha=5\%$ and we study the empirical power function.
In the second part, we present an application to some economical
series.

\subsection{Methodology}
On the contrary to the estimation problem, a smooth wavelet is not needed.
The test statistic is then computed with the Haar wavelet since it
has a small support and then it leads to a fast computation time.
The critical value of the test is determined with bootstrap methods:
the standard deviation of the test statistic is computed thanks to
$N_{boot}=20$ resampling.  The size of the simulated samples is
$n=2048$ which is reasonable for bi-dimensional problems in an
asymptotic context. For the real life data example, the number of
data is around $n=4000$. For the simulation part, the empirical
level of the test is derived from $N_{MC}=500$ replications for each
test problem.

\subsection{Simulations}
 The setup of our simulations is closely related to the work of
 \citep{Genest/Remillard/Beaudoin:2008},
 except that they consider small samples (of size $150$) since their  test procedures
 are based on the empirical
copula distribution (and thus generate parametrical rates). To
explore various degrees of dependance, three values of Kendal's tau
are considered, namely $\tau = 0.25, 0.50, 0.75$ for the following
copula families: Clayton, Gumbel, Frank, Normal and Student with
four degrees of freedom (df). Calculations are made with  the MatLab
Sofware.  The results of the simulations are presented in Table
\ref{outside}. For an easier reading,  the estimated standard errors
of the empirical powers are presented in italics. Furthermore, for
each testing problem we highlighted the estimated errors of the
first type (estimators of $\alpha=0.05$) using bold characters. In
brackets, we give the results obtained by
\citep{Genest/Remillard/Beaudoin:2008} with  their test
procedures, denoted CvM and built on rank-based versions of the
familiar Cram\'er-von Mises statistics. It would be also possible,
if one is interested in,  to compare with  the different  test procedures
(based on the empirical copula distribution) proposed
also by \citep{Genest/Remillard/Beaudoin:2008}.

Let us now summarize the conclusions made from the simulation results.
\begin{itemize}
\item Our test is  degenerated: we almost
always accept $H_0$ (when $H_0$ is true) while the procedure of
Genest et al. \citep{Genest/Remillard/Beaudoin:2008} produces an
excellent estimation of the prescribed level $\alpha$. It is a
characteristic of the adaptive minimax procedures.
\item For small level of dependence $\tau=0.25$, our procedure is very
competitive and produces (almost) always a better empirical power
than the CvM test. The results are spectacular when the fit
$c_{\lambda_0}$ is a Student(4).
\item
When a large Kendal's tau is considered, our procedure fails when
the data are issued from a Clayton copula density. The procedure is
not available to recognize a structure of dependence modeled with a
Clayton.
\item The improvment of our results with respect to the CvM test is
decreasing with the Kendal's tau. The CvM test becomes better when
the tau is increasing whereas  for us it is the opposite.
\end{itemize}
In conclusion, we recommend the use of our test procedures when the Kendal's tau is
not too large since it seems to outperform the existing procedures
based on the copula distribution. This situation  corresponds to our
theoretical setup related to the functional spaces in which  the
unknown copula density is supposed to live. Unfortunately, the practical
results do not give hope for using this procedure when the copula
densities present high peaks (as it is the case for the Clayton copula density with a
large tau).

\subsection{Real data}
We present now an application  to real data of our test procedure.
The level of each test (with simple null hypothesis or multivariate
null hypothesis) is $\alpha=5\%$. To obtain the empirical level,
$N=50$  replications of our procedure computed with the half of the
available data (chosen randomly) is used.
 Table \ref{test_frees}  gives the empirical probability to reject
 the null hypothesis and the final decision. "Yes" means that we accept that the structure of
 dependence belongs to the considered family and "No" that we reject
 the fitting.

We consider the data of \citep{Frees/Valdez:1998}, which were also
analyzed by \citep{Genest/Ghoudi/Rivest:1998},
 \citep{Klugman/Parsa:1999}, \citep{Chen/Fan:2005} and
\citep{Genest/Quessy/Remillard:2006}, among others. The data consist
of the indemnity payment (LOSS) and the allocated loss adjustment
expense (ALAE) for 1466 general liability claims.

 We consider the following test problems:
$$H_0: c\in {\cal{C}}_\Lambda$$ where the parametrical family ${\cal{C}}_\Lambda$
is described in Table \ref{test_frees}. Since the Kendall's tau
 computed with the sample is $\tau=0.31$, we choose an adapted grid of parameters
 for each parametrical family of copula densities.
 Next, assuming that the density copula of the data belongs to a
fixed parametric family, we estimate the parameter $ \lambda$
\begin{itemize}
\item by $\hat \lambda$ in inverting  the Kendall's tau (third part of Table
\ref{test_frees} where $H_0: c=c_{\widehat{\lambda}}$).
\item by $\tilde \lambda$ in minimizing the average square error (ASE)
 computed thanks to the benchmark given in Figure \ref{fig3}
 (fourth part of Table \ref{test_frees} where $H_0: c=c_{\widetilde {\lambda}}$). For information, we give the
relative $ASE$ computed with $c_{\widetilde {\lambda}}$ into
brackets.
\end{itemize}
The various authors who analyzed this data set concluded that the
Gumbel copula provides an adequate representation of the underlying
dependence structure.  The Gumbel parametric family of extreme-value
copulas captures the fact that almost all large indemnity payments
generate important adjustment expenses (e.g., investigation and
legal costs) while the effort invested in the treatment of a small
claim is more variable. Accordingly, the copula exhibits positive
but asymmetric dependence. Confirming this result, the adaptive
method of estimation proposed by
\citep{Autin/LePennec/Tribouley:2008} provides a benchmark
 (see Figure \ref{fig3}) for the copula density
associated with the data.

\section{Discussion}

The paper is mainly devoted to construct an optimal  procedure for
solving a general nonparametric problem of test: both hypotheses are
composite, very general parametric family could be considered under
the null. Our procedure is proved asymptotically to be adaptive
minimax and the minimax separating rate is exhibited over a range of
Besov balls.

Thanks to the simulations and a application to real data,
our procedure seems to be competitive on the power point of view
even if the setting of test under consideration is, in the simulation study, clearly
parametric.

It is worthwhile to point out that the copula model requires more
regularity (than the usual density model) since the approximation
due to the rank-based statistics needs to be accurate enough (see
Lemma \ref{nombredei}).

One must notice that only copulas densities belonging to {\bf dense
Besov spaces} (i.e. defined with a parameter $p$ larger than 2) are
under consideration in this paper although
 several copula densities with a strong
dependence structure belong to  {\bf sparse Besov spaces} (i.e.
defined with a parameter $p$ smaller than 2). As it is illustrated
in the simulation study, our test procedure fails for the Clayton
copula density with large parameters. This  density is suspected to
belong to a sparse Besov ball. The study  of  sparse Besov balls
would require the determination of a new test strategy which would
lead to another minimax rate of testing: these objectives are beyond
those of the present paper and will be explored in  a further work
since the set of  copulas densities contains a number of sparse
functions. For sparse Besov balls and in the white noise model for
testing the existence of the signal, \citep{Lepskii/Spokoiny:1999}
proved that the minimax  testing rate in the sparse and the dense
cases
 is different. They also proved that it is possible to built an adaptive minimax (non
linear) procedure of test for the sparse case.

A very close problem is the sample comparison test
(problem with two samples). It could be interesting to test if the
structure of dependence between a couple of variables $V_1=(X,Y)$ is
the same as for another couple $V_2=(Z,T)$.   This problem of tests could be stated as follows:

\begin{eqnarray*}
H_0 : c_{V_1}=c_{V_2} \quad &\mbox{ {\rm against} }&
 \quad  H_1:(c_{V_1},c_{V_2}) \in \Gamma(v_n(\tau)),
\end{eqnarray*}
with
\begin{eqnarray*}
\Gamma(v_n)=\left \{ c_{V_1}\in b_{s_1,p_1,\infty}(M_1)\right\} \cap
\left\{c_{V_2}\in
b_{s_2,p_2,\infty}(M_2)\right\}\\
\cap \left\{(c_{V_1},c_{V_2}): \|c_{V_1}-c_{V_2}\| \geq v_n.
\right\}\end{eqnarray*} where $v_n$ is the separating rate of both
hypotheses. In an analogous way as in Section 3,  the rule for the comparison test would be
\begin{eqnarray*}
 D&=&  \I_{ \displaystyle{ \max_{j\in J}
(\sum_k \widetilde {\theta_{j,k}}- t_j  )}> 0}
 \end{eqnarray*}
with
 \begin{eqnarray*}
 \widetilde {\theta_{j,k}}
 &=&
 \frac{1}{n_2(n_2-1)}\mathop{\sum_{i_1,
i_2\in {\cal I}_2}}_{ i_1 \neq i_2 }
\left(\phi_{j,k}\left(\frac{R^X_{i_1}}{n_1} ,\frac{R^Y_{i_1}}{n_1}
\right)- \phi_{j,k}\left(\frac{R^Z_{i_1}}{n_1}
,\frac{R^T_{i_1}}{n_1} \right) \right)\\&&\hspace{3cm} \times
\left(\phi_{j,k}\left(\frac{R^X_{i_2}}{n_1}, \frac{R^Y_{i_2}}{n_1}
\right)- \phi_{j,k}\left(\frac{R^Z_{i_2}}{n_1},
\frac{R^T_{i_2}}{n_1} \right)\right),
\end{eqnarray*}
where $R^X,R^Y,R^Z,R^T$ are the rank statistics associated with
$X,Y,Z,T$. Using the same tools as in
\citep{Butucea/Tribouley:2006}, in which  the homogeneity in law of
the both samples is studied, it is possible to prove that this test
is adaptive optimal and that the minimax rate of testing is
$$
v_n=\left(\frac{n}{\sqrt{\log(\log(n_2))}}\right)^{-2(s_1\wedge
s_2)/(4(s_1\wedge s_2)+2)}.
$$
Obviously, all these test procedures could be used in the
multivariate framework ($d >2$),  but as usual in the nonparametric context,
it will  provide slower minimax rates of testing.

\section{Proof of Theorem \ref{BSup}}
Recall that for any given $\lambda\in \Lambda$,  $\P_{\lambda}$
(respectively $\P_{c}$) denote the distribution associated with
density $c_\lambda$, respectively with $c$. In the same spirit, denote
also $\E_{\lambda}$ and $\V_{\lambda}$ (respectively $\E_{c}$ and
$\V_{c}$) the expectation and the variance with respect to
$\P_{\lambda}$, respectively to
 $\P_{c}$. When no index appears in $\E$ or in $\P$ it means that the underlying distribution is either $\P_{c}$ or
 $\P_{\lambda}$.

\subsection{Expansion of the statistics of interest}
Fix a level $j$ in $J$. For the test problem, the statistic of
interest $\widetilde{ T_j}(\lambda)$ (for $\lambda\in\Lambda$)
defined in (\ref{Ttilde}) is an estimate of
\begin{eqnarray*}\label{expansionVRAIE}
T_j(\lambda) &=&\sum_k \theta_{j,k}(\lambda)
=\sum_k\left(c_{j,k}-c_{j,k}(\lambda)\right)^2,
\end{eqnarray*}
which is the quantity that we need to detect under the alternative.
 It would be  useful to expand the statistic
$\widetilde{T_j}(\lambda)$ as follows
\begin{eqnarray}
&&\widetilde{T_j}(\lambda)= 2 T_j^\diamond (\lambda)
+{T_j}^\heartsuit +T_j^\spadesuit +2 T_j^\clubsuit
(\lambda)+T_j(\lambda) \label{expansion}
\\&&\quad \quad  =2\sum_k {\theta}_{j,k}^\diamond(\lambda)+\sum_k
\theta_{j,k}^\heartsuit+\sum_k {\theta}_{j,k}^\spadesuit +2\sum_k
{\theta}_{j,k}^\clubsuit (\lambda)+\sum_k \theta_{j,k}(\lambda),
\nonumber
\end{eqnarray}
where
\begin{eqnarray*}
{\theta}_{j,k}^\heartsuit & =&
\frac{1}{n_2(n_2-1)}\mathop{\sum_{i_1, i_2 \in {\cal I}_2}}_{ i_1
\neq i_2 }\left(\phi_{j,k}(F(X_{i_1}), G(Y_{i_1}))-c_{j,k}\right)\\
&&\times (\phi_{j,k}\left( F(X_{i_2}), G(Y_{i_2}))-c_{j,k}\right)
\\
\theta_{j,k}^\spadesuit
 &=&\frac{1}{n_2(n_2-1)}\mathop{\sum_{i_1,
i_2 \in {\cal I}_2}}_{ i_1 \neq i_2} \left(\phi_{j,k}\left(\frac{R_{i_1}}{n_1},
\frac{S_{i_1}}{n_1} \right)-
\phi_{j,k}\left(F(X_{i_1}),G(Y_{i_1})\right)
\right)\,\\&&\times\left(\phi_{j,k}\left(\frac{R_{i_2}}{n_1},
\frac{S_{i_2}}{n_1}
\right)-\phi_{j,k}\left(F(X_{i_2}),G(Y_{i_2})\right) \right)
\\
\theta_{j,k}^\clubsuit (\lambda)
 &=&\frac{1}{n_2(n_2-1)}\mathop{\sum_{i_1,
i_2 \in {\cal I}_2}}_{ i_1 \neq i_2
}\left(\phi_{j,k}\left(\frac{R_{i_1}}{n_1}, \frac{S_{i_1}}{n_1}
\right)-\phi_{j,k}\left(F(X_{i_1}),G(Y_{i_1})\right)
\right)\\&&\times\left(\phi_{j,k}\left(F(X_{i_2}),G(Y_{i_2})
\right)- c_{j,k}(\lambda) \right)\\
{\theta}_{j,k}^\diamond (\lambda) &=& \frac{1}{n_2}\sum_{i_1\in
{\cal I}_2} \left(\phi_{j,k}(F(X_{i_1}), G(Y_{i_1}))-c_{j,k}\right)
(c_{j,k}-c_{j,k}(\lambda)).
\end{eqnarray*}
The sequence $\{c_{j,k}\}_{j,k}$ denotes the unknown scaling
coefficients of the unknown copula density $c$. Recall that
$$
\widehat{T_j}(\lambda)=\sum_k \widehat{\theta_{j,k}}(\lambda),
$$
with
\begin{eqnarray*}
\widehat{\theta_{j,k}}(\lambda) &=&
\frac{1}{n_2(n_2-1)}\mathop{\sum_{i_1, i_2 \in {\cal I}_2}}_{ i_1
\neq i_2 }\left(\phi_{j,k}(F(X_{i_1}),
G(Y_{i_1}))-c_{j,k}(\lambda)\right) \\
&&\times \left(\phi_{j,k}( F(X_{i_2}),
G(Y_{i_2}))-c_{j,k}(\lambda)\right).
\end{eqnarray*}
 The following lemma gives some evaluation
for the first moments of each statistic of interest.

\vspace{0.2cm}

\begin{lemma}\label{etude-stat}
Let $q$ be a positive integer and assume that $\phi$ is continuously
$q-$differentiable. Let $j$ be a level smaller than $ j_\infty$
defined in (\ref{j0jinfty}). Then, it exists some positive
constant $\kappa$ which may depend on $\phi$, $\|c\|_\infty$,
$\|c_\lambda\|_{\infty}$ and $M$ such that
\begin{eqnarray*}
\E \widehat{T_{j}}(\lambda)&=&T_j(\lambda)\quad\mbox{ {\rm and} }\quad \V
\widehat{T_{j}}(\lambda)\leq
\kappa  \left(\left(\frac{2^{j}}{n_2}\right)^2+
\left(\frac{2^{j}}{n_2}\right)\, T_j(\lambda)\right) \\
\E |T_{j}^\spadesuit |&\leq& \kappa\;\frac{\log(n_1)}{n_1}\\
 \E_c
|T_{j}^\clubsuit (\lambda)|&\leq& \kappa\,
\left(\frac{\log(n_1)}{n_1}T_j(\lambda)\right)^{1/2}\;\mbox{ {\rm
and}  }\; \E_\lambda( T_{j}^\clubsuit  (\lambda))^2\leq \kappa
\,2^{j}\left(\frac{\log(n_1)}{n_2n_1}\right).
\end{eqnarray*}
\end{lemma}
Using the Bernstein Inequality, we establish the following bound for
the deviation of the statistic $T_{j}^\diamond(\lambda)$ under the alternative. The proof
is postponed to Appendix B.

 \vspace{0.2cm}

\begin{lemma}\label{Bernstein}
  For any level $j$, for all $x>0$
\begin{eqnarray*}
\P_c\left(\,|T_{j}^\diamond(\lambda)|\,\geq x\right)\leq
\exp\left(-K (\frac{n_2^2x^2}{n_2 T_j(\lambda)+n_2x2^{j}T_j(\lambda)^{1/2}})\right),
\end{eqnarray*}
  where  $K$ is a positive
constant depending on $L,\|\phi\|_\infty $ and
$\|c\|_\infty$.
\end{lemma}
Using a result from \citep{Gine/Latala/Zinn:2000}, we establish the
following bound for the deviation of the $U$-statistics
$\widehat{T_j}(\lambda)$ and $T_j^\heartsuit$. The proof is
postponed to Appendix C.

\vspace{0.2cm}

\begin{lemma}\label{LDUstat}
For any level $j$,
as soon as $ x\geq 2^j\,n_2^{-1}\sqrt{\log(\log(n_2))}$, for all
$\mu>0$,
\begin{eqnarray*}
\P_\lambda\left(|\widehat{T_j}(\lambda)|> \mu x \right)+\P_c\left(|
T_j^\heartsuit|> \mu x \right)&\leq& K_g(\log(n_2))^{-\delta}
\end{eqnarray*}
for any positive
 $\delta \leq \mu^2(K_g K_1)^{-1}$,
  where $K_g$ is an universal positive
constant given in \citep{Gine/Latala/Zinn:2000} and $K_1$ is a
positive constant depending on $L,\|\phi\|_\infty $ and either
$\|c_{\lambda}\|_\infty$ or  $\|c\|_\infty$ depending on  the
underlying distribution i.e. either $\P_\lambda$ or $\P_c$.
\end{lemma}

\subsection{Proof of Relation (\ref{PremE}) (First type error)} Let us fix
$\lambda\in\Lambda$ and set $$ p_\lambda = \P_{\lambda}
\left(\max_{j \in J} \left[\inf_{\lambda' \in \Lambda}
\widetilde{T_j}(\lambda') - t_j\right]
>0 \right).$$
Notice that under the null
$$
 T_j^\diamond (\lambda)
 =T_j(\lambda)=0 \mbox{ and }
 {T_j}^\heartsuit=\widehat{T_j}(\lambda).
$$
Using  expansion (\ref{expansion}), we get
\begin{eqnarray*}
p_{\lambda}& \leq &   \sum_{j \in J} \P_{\lambda}\left(
\inf_{\lambda^\prime \in \Lambda} \widetilde{T_j}(\lambda^\prime) > t_j\right)\nonumber  \\
& \leq &   \sum_{j \in J} \P_{\lambda}\left(
 \widetilde{T_j}(\lambda) > t_j\right)\nonumber  \\
& \leq &   \sum_{j \in J}
\left\{\P_{\lambda}\left(|\widehat{T_j}(\lambda) |>
\frac{t_j}{3}\right) + \P_{\lambda} \left(|T_j^\spadesuit |
> \frac{t_j}{3}\right) + \P_{\lambda} \left(|T_j^\clubsuit (\lambda)| >
\frac{t_j}{3}\right)\right\} \nonumber
\end{eqnarray*}
Due to Lemma  \ref{etude-stat} and using Markov Inequality, we obtain
\begin{eqnarray*}
p_{\lambda}& \leq & \sum_{j\in J}
\P_{\lambda}\left(|\widehat{T_j}(\lambda) |> \frac{t_j}{3} \right) +
\sum_{j\in J} \left\{\frac{\E_{\lambda}| T_j^\spadesuit |}{
(t_j/3)}+\frac{ \E_{\lambda}( T_j^\clubsuit (\lambda))^2}{
(t_j/3)^2}\right\} \nonumber \\
&\leq& \sum_{j\in J} \P_{\lambda}\left(|\widehat{T_j}(\lambda) |>
\frac{t_j}{3}\right)\\&& + K\sum_{j\in
J}\left\{\left(t_j/3\right)^{-1}\frac{\log(n_1)}{n_1}+\left(t_j/3\right)^{-2}\left(
\frac{2^{j} \log(n_1)}{ n_1n_2 }\right)\right\}.
 \label{THM1firstIneq}
\end{eqnarray*}
Notice that $\widehat{T_j}(\lambda)$ is centered under $\P_\lambda$,
then applying Lemma \ref{LDUstat}, where $t_j$  is $t_j= 3\mu
\;2^j\,n_2^{-1}\sqrt{\log(\log(n_2))}$,  the constant  $\mu$ is
defined in (\ref{seuil}) and since $\mbox{{\rm card}} (J) \leq \log
(n_2)$, one obtains
\begin{eqnarray*}
p_{\lambda}& \leq & K_g\mbox{{\rm card}} (J)\; (\log(n_2))^{-\delta} + K \mbox{{\rm card}} (J)
2^{-j_0}\left(\frac{ n_2\log(n_1)}{ n_1\sqrt{\log\log(n_2)} }\right)\\
& & +K \mbox{{\rm card}} (J)
2^{-j_0}\left( \frac{\log(n_1) n_2^2}{ n_2 n_1 \sqrt{\log\log(n_2)}}\right)\\
& \leq &    K_g  (\log(n_2))^{1-\delta}   + K  2^{-j_0} \left(\frac{\log(n_1)\log(n_2)}{\sqrt{\log\log(n_2)} }\right),
\end{eqnarray*}
where the last inequality holds since $\delta$ satisfies $\delta  \leq \mu^2(2 K_g K_1)^{-1}$
(see Lemma \ref{LDUstat}).  Since $\mu$ is such that  $\mu >
\sqrt{2 K_g K_1}$, relation (\ref{PremE}) is proved if one takes $\delta=\mu^2(2 K_g K_1)^{-1}$.

\subsection{Proof of Relation (\ref{SecondE}) (Second type error)}
Let us fix $\tau \in {\cal S}$ and
$c\in\Gamma(Av_{nt_n^{-1}}(\tau))$ and set
$$
 p_c=\P_{c} \left(\max_{j \in J} \inf_{\lambda \in
\Lambda} \widetilde{T_j}(\lambda) - t_j \leq 0 \right).
$$
Using the expansion (\ref{expansion}), we get, for any $j^\star\in
J$
\begin{eqnarray}
p_{c}& \leq & \P_c \left( \inf_\lambda \left\{2 T_{j^\star}^\diamond
(\lambda)+T_j(\lambda)+ T_{j^\star}^\heartsuit +
T_{j^\star}^\spadesuit +2T_{j^\star}^\clubsuit (\lambda)\right\}\leq
t_{j^\star}
\right)\nonumber\\
 & \leq & \P_c \left( \inf_\lambda \left\{2
T_{j^\star}^\diamond (\lambda)+T_{j^\star}(\lambda)\right\}\leq 2t_{j^\star}
\right)\nonumber\\
&&+ \P_c \left(T_{j^\star}^\heartsuit +  T_{j^\star}^\spadesuit +2
\inf_\lambda \left\{T_{j^\star}^\clubsuit (\lambda)\right\}\geq
t_{j^\star}\right)\nonumber\\
& \leq & \P_c \left( \inf_\lambda \left\{2 T_{j^\star}^\diamond
(\lambda)+T_{j^\star}(\lambda)\right\}\leq 2t_{j^\star} \right)+ \P_c
\left(T_{j^\star}^\heartsuit \geq t_{j^\star}/3\right)\nonumber\\
&&
 + \P_c \left(T_{j^\star}^\spadesuit \geq t_{j^\star}/3\right)
+\P_c \left(  \inf_\lambda \left\{T_{j^\star}^\clubsuit
(\lambda)\right\}\geq t_{j^\star}/6\right)
\nonumber\\
&=&p_{c1}(j^\star)+p_{c2}(j^\star)+p_{c3}(j^\star)+p_{c4}(j^\star).\label{2proba}
\end{eqnarray}
Let us explain how $j^\star$ is chosen. From the wavelet expansion
(\ref{decomposition}) and Lemma \ref{etude-stat}, one has
$$\E_c\widehat{T_{j^\star}}(\lambda)=T_{j^\star}(\lambda)=\int(c-c_\lambda)^2-B_{j^\star}(\lambda),$$
 where  $T_{j^\star}$, $B_{j^\star}$ are
defined in (\ref{TB})  and $t_j^\star$ is the critical value  given
in (\ref{seuil}). Since $c$ is in $\Gamma(Av_{nt_n^{-1}}(\tau))$ and
 $c_\lambda$ lies in $b_{s_\Lambda, p_\Lambda , \infty}(M_\Lambda)
 \subset b_{s, p , \infty} (M)$, the function $(c-c_\lambda)$ is in
$b_{s,p,\infty}(M)$. We can choose $j^{\star}$ such that
\begin{eqnarray*}\label{jstar}
2^{j^\star} & =& \left(\frac{{\tilde
K}}{3\mu}\frac{n_2}{\sqrt{\log\log(n_2)}}\right)^{1/(2s+1)},
\end{eqnarray*}
which is possible due to our choice of  $j_\infty$ and because
$s\geq 1/2$; the constant ${\tilde
K}$ appears in (\ref{reg}). It implies that $B_{j^\star}\leq t_{j^\star}$ since
$B_{j^{\star}} \leq {\tilde K} 2^{-2j^{\star}s}$ (see Inequality
(\ref{reg})). Next, since $c\in \Gamma(Av_{nt_n^{-1}}(\tau))$, one
has $\int(c-c_{\lambda^\prime})^2\geq A^2(v_{nt_n^{-1}}(\tau))^2$
for all $\lambda^\prime  \in \Lambda$. Focusing on rates
$v_{nt_n^{-1}}(\tau)$ combined with positive constant $A$ which
satisfy $4t_{j^\star}\leq (A v_{nt_n^{-1}}(\tau))^2$, one obtains
\begin{eqnarray}\label{tT}
\frac{t_{j^\star}}{\E_c\widehat{T_{j^\star}}(\lambda)}\,=\,\frac{t_{j^\star}}{T_{j^\star}(\lambda)}&\leq
&1/3.
\end{eqnarray}
Coming back  to the evaluation of the probability terms (see relation (\ref{2proba})),  we first consider $p_{c1}(j^\star)$.
Consider an $\eta$-net $\Lambda_{\eta}$  on the set $\Lambda$ that is  for any $\lambda$ in $\Lambda$, denote ${\tilde \lambda}$ the closest (in the Euclidean sense) element in $\Lambda_{\eta}$ to $\lambda$ (closer than $\eta$). Due to assumption {\bf A0}, let us prove that for any $j \in J$,
$T_{j}^\diamond
({\tilde \lambda})+T_{j}({\tilde \lambda})$ is close to
 $T_{j}^\diamond
(\lambda)+T_{j}(\lambda)$:
\begin{eqnarray}
|T_{j}^\diamond ({\tilde \lambda})-T_{j}^\diamond ({\lambda})| & = &
\left | \sum_k \left [ \frac{1}{n_2} \sum_{i_1 \in {\cal I}_2}
\left(\phi_{j,k}(F(X_{i_1}), G(Y_{i_1}))-c_{j,k} \right)
\left(c_{j,k}({\tilde \lambda})-c_{j,k}(\lambda) \right)  \right ]\, \right|\nonumber \\
&\leq & \left[ \sum_k \left(c_{j,k}({\tilde
\lambda})-c_{j,k}(\lambda)\right)^2 \right]^{1/2}\nonumber \\
&&\hspace{1cm}\times \,\left[ \sum_k \left(\frac{1}{n_2}
\sum_{i_1 \in {\cal I}_2} \left(\phi_{j,k}(F(X_{i_1}), G(Y_{i_1}))-c_{j,k}\right) \right)^2\right]^{1/2}\nonumber \\
&\leq & Q \eta^{\nu} \; 2^{2j}(2\|\phi\|_{\infty} + \|c\|_{\infty}
2^{-4j})^{1/2}. \nonumber
\end{eqnarray}
In the same way, one has,
\begin{eqnarray*}
|T_j(\lambda) - T_j({\tilde \lambda})| & \leq &
 \kappa Q \eta^{\nu}\;\left(2\mbox{{\rm max}}(\|c\|, \|c_{{\tilde \lambda}}\|) +  Q \eta^{\nu}\right).
\end{eqnarray*}
Choosing $\eta=n^{-b}$ with $b \nu >1$, then by (\ref{tT}) and
applying Lemma \ref{Bernstein}, we get
\begin{eqnarray}
p_{c1}(j^\star)&\leq& \sum_{\lambda\in\Lambda_{n^{-b}}} \P_c\left( 2
T_{j^\star}^\diamond
(\lambda)+T_{j^\star}(\lambda)\leq 2t_{j^\star}\right)\nonumber\\
&\leq& \sum_{\lambda\in\Lambda_{n^{-b}}} \P_c\left(  T_{j^\star}^\diamond
(\lambda)\leq -T_{j^\star}(\lambda)/6\right)\nonumber\\
&\leq&
\sum_{\lambda\in\Lambda_{n^{-b}}}\exp\left[-K\,\left(\frac{n_2}{2^{2j^\star}}T_{j^\star}(\lambda)\wedge
\frac{n_2}{2^{j^\star}}T_{j^\star}(\lambda)^{1/2} \right)\right]\nonumber\\
&\leq&
\sum_{\lambda\in\Lambda_{n^{-b}}}\exp\left[-K_2\left(n_2 t_{j^\star}\wedge
\frac{n_2}{2^{j^\star}}t_{j^\star}^{1/2} \right)\right]\nonumber\\
&\leq &(\frac{D(\Lambda)}{n^{-b}})^{d_{\Lambda}}\; \exp\left[
-K_3\,\left(2^{j^{\star}}(\log\log(n_2))^{1/2}\right.\right.\nonumber\\
&&\hspace{4cm}   \left.\left. \wedge \left(\frac{ n_2}
{2^{j^{\star}/2}}\right)^{1/2}(\log\log(n_2))^{1/4}\right)
\right]\label{2ndeProbaT3}
\end{eqnarray}
where $D(\Lambda)$ is the diameter of $\Lambda$ and $K$, $K_2$ and
$K_3$ are positive constants. Both terms behind the minus sign in
the exponential of the right hand side of the last inequality tend
to infinity with a power of $n$ since $s > 1/2$.
This implies that $p_{c1}(j^\star)$ goes to zero  as $n$ goes to infinity.

 \vspace{0.5cm}
Now, it remains  to verify that $p_{c2}(j^\star),p_{c3}(j^\star)$
and $p_{c4 }(j^\star)$ are going to zero as $n_1\wedge n_2$ goes to
infinity. Using again the bound (\ref{tT}), Lemma \ref{LDUstat} for
some positive $\delta$, Lemma \ref{etude-stat} and the definition of
the critical value (\ref{seuil}), one gets

$p_{c2}(j^\star)+p_{c3}(j^\star)+p_{c4}(j^\star)$

\begin{eqnarray}
&&  \leq  \P_c \left(T_{j^\star}^\heartsuit \geq
t_{j^\star}/3\right)+
 + \P_c \left(T_{j^\star}^\spadesuit \geq t_{j^\star}/3\right)+
 \P_{c}\left( T_{j^\star}^\clubsuit
(\lambda)\geq t_{j^\star}/6\right)\nonumber\\
&& \leq K_g (\log (n_2))^{-\delta}+ 9
\,\frac{\E_c|T_{j^\star}^\spadesuit (\lambda)|^2 }{t_{j^\star}^2}+
\,6\,\frac{\E_c|T_{j^\star}^\clubsuit |
}{t_{j^\star}} \nonumber \\
&& \leq K_g (\log (n_2))^{-\delta}+ 9 \kappa
\frac{2^{j^\star}\log(n_1)}{n_1n_2\,t_{j^\star}^2}
+\,6 \kappa \frac{\log(n_1)}{n_1\,t_{j^\star}}\nonumber\\
&&\quad \leq K_g (\log (n_2))^{-\delta}+ 6\kappa \;
\left(2^{-j^\star}\frac{n_2}{n_1}\frac{\log(n_1)}{\sqrt{\log\log(n_2)}}\right)\label{2ndeProbaT2},
\end{eqnarray}
which tends to zero with our choice of $j^\star$
and where  $\kappa$ is the positive constant appearing in Lemma
\ref{etude-stat}. Inequalities (\ref{2ndeProbaT3}) and
(\ref{2ndeProbaT2})  entail that the right hand side of
(\ref{2proba}) is less than any $\alpha \in (0,1)$ as $n$ is large
enough. To finish the proof, observe that the choice of
$v_{nt_n^{-1}}(\tau)$ is driven by the fact that it corresponds to
the smallest sequence such that $4 t_j^{\star} \leq (A
v_{nt_n^{-1}}(\tau))^2$, which leads to
$$v_{nt_n^{-1}}(\tau)\geq \left(2^{j^\star} \mu \frac{\sqrt{\log(\log n_2)}}
{n_2 }\right)^{1/2} \geq \left(\frac{n_2}
{\sqrt{\log\log(n_2)}}\right)^{2s/(4s+2)}.$$


\section{Proof of Theorem \ref{BInf}}
Without loss of generality, we suppose that the support of the
scaling function $\phi$ and its associated wavelet function $\psi$
is $[0,1]$. Moreover recall that $\int_0^1 \psi^{\epsilon} =0$. Let
us give some $a>0$ which must be small enough.

\subsection{ Discretisation of ${\cal S}$}
For any given $\tau=(s,p,M) \in {\cal S}$, denote by $j(\tau)$ the level
$$
2^{j(\tau)}=(nt_n^{-1})^{2/(4s+2)}
$$
and define $s_j$  the solution of the equation $j=j(s_j,p,M)$ for
any resolution level $j \in {\tilde J}=\{j_{s_{\mbox{max}}},\ldots, j_{s_{\mbox{min}}}\} \subset \{j_0,\ldots,j_{\infty}\}$ with
 $$j_{s_{\mbox{max}}}=\lfloor j(s_{\mbox{max}},p,M)\rfloor \mbox{ and } j_{s_{\mbox{min}}}=\lfloor
j(s_{\mbox{min}},p,M)\rfloor.
$$

Consider now  the set ${\cal S}_n=\{\tau_j=(s_j,p,M), j\in {\tilde J}\}$
which appears as a discretisation version of a subset of ${\cal S}$ whose
cardinality is of order $O(\log (n))$.

\subsection{Prior and parametric family included in the alternatives}
For any $s_j \in {\cal S}_n$, define a prior $\pi_j$ which is concentrated on the class of the random functions
$$c_j(u,v) = c_{\lambda_0} (u,v) +  \sum_k  \sum_{\epsilon =1}^3 \delta_{k} u_{j}(n) \psi^{\epsilon}_{j,k}(u,v),$$
where $c_{\lambda_0}$ is defined in assumption {\bf AInf} and
$$
P(\delta_{k} =1 )=P(\delta_{k} =-1 )=1/2\quad \mbox{ and }\quad
u_{j}(n)= C_1 M (nt_n^{-1})^{-\frac{2(s_j+1)}{4s_j+2}}
$$
for $C_1$  such that $3 M^2 C_1^2=2a^2$. Let $j$ be any index in ${\tilde J}$.
Since $\int \psi=0$ and when $a$ is small enough (to guarantee that
$c_j\geq 0$), $c_j$ is a density. Easy calculations imply that
 \begin{eqnarray*}
\| c_j-c_{\lambda_0}\|^2 = M^2 C_1^2   \; (v_{nt_n^{-1}}(\tau_j))^2
> a^2 \;(v_{nt_n^{-1}})^2\label{L2} .\end{eqnarray*} Moreover,
if $a$ is small enough, we have $3\; C_1^p <1$ and
\begin{eqnarray*}
2^{j(s_j+1-2/p)p} \sum_k \sum_{\epsilon} | \int
c_j\psi_{jk}^\epsilon |^p&=&
 2^{j(s_j+1-2/p)p} \sum_k \sum_{\epsilon} | u_{j}(n) |^p\\
 & =&  3  C_1^p M^p \leq  M^p,
\end{eqnarray*}
implying that $c_j\in b_{s_j,p,\infty}(M)$. Denote by  ${\cal A}_{j,n}
(a)$ the set of densities
$${\cal A}_{j,n}
(a)=
\{c \in b_{s_j,p,\infty}(M) : \;
\displaystyle{\inf_{\lambda \in \Lambda}}
 \|c - c_{\lambda}\|^2 > a^2 (v_{nt_n^{-1}}(\tau_j))^2 \}.
 $$
and consider the variation between both distributions
$\P_{\lambda_0}$ and $\P_\Pi$
$$
Var(\P_{\lambda_0}, \P_\Pi)= \frac 1 2 \int
\left|\frac{d\P_{\Pi}}{d\P_{\lambda_0}}-1\right| d\P_{\lambda_0},
$$
where
$$
\frac{d\P_{\Pi}}{d\P_{\lambda_0}} = \frac{1}{N_n } \sum_{j \in {\tilde J}}
 \frac{d\P_{j}}{d\P_{\lambda_0}}=\frac{1}{N_n } \sum_{j \in {\tilde J}}
  \E_{\pi_j}^{(n)} [\frac{c_{j}}{c_{\lambda_0}}],
  $$
  and $N_n=\mbox{ {\rm card} }({\tilde J})$.
Assuming that the following assertion holds
\begin{eqnarray}
\lim_{n \rightarrow \infty} \inf_{j \in {\tilde J}} \pi_j (c \in {\cal
A}_{j,n} (a))=1 ,\label{alternative}
\end{eqnarray}
we deduce that the left hand side ($LHS$) of relation (\ref{BI}) without
the limit is bounded from below by
\begin{eqnarray*}
LHS&\geq& \P_{\lambda_0}( D_n=1) + \sup_{\tau_j \in {\cal S}_n} \sup_{c
\in {\cal A}_{j,n} (a)} \P_c ( { D}_n =0)\nonumber\\
 & \geq & 1 - Var(\P_{\lambda_0}, \P_\Pi)(1 + o_n(1)),\label{Variation}
\end{eqnarray*}
as $n$ large enough.  Since the supports of the functions $c_j$ and
$c_{j'}$ are disjoint for $j \neq j'$, one has
\begin{eqnarray*}
1 - Var(\P_{\lambda_0}, \P_\Pi) & \geq & 1 - \frac 1 2
\frac{1}{N_n^2} \sum_{j \in {\tilde J}}\E_{\lambda_0}
 \left[ \left(\int \prod_{i=1}^n \frac{c_j(U_i,V_i)}{c_{\lambda_0}(U_i,V_i)} d\pi_j(c_j)\right)^2 -1 \right]  \nonumber\\
&\geq& 1 -o_n(1)
 \end{eqnarray*}
provided that
\begin{eqnarray}
 \lim_{n \rightarrow \infty} \frac{1}{N_n^2} \sum_{j \in {\tilde J}}
  \E_{\lambda_0} \left[ \left(\int \prod_{i=1}^n \frac{c_j(U_i,V_i)}{c_{\lambda_0}
  (U_i,V_i)} d\pi_j(c_j)\right)^2 \right] = 0. \label{vrais}
\end{eqnarray}
Relation (\ref{BI}) is thus proved if (\ref{alternative}) and
(\ref{vrais}) are satisfied.  The remaining proofs are given in the
sequel.

\subsection{Proof of Relation (\ref{alternative})}
Let $\Lambda^\prime$ be a subsect of $\Lambda$. We have
\begin{eqnarray*}
\pi_j\left(\inf_{\lambda \in \Lambda} \|c_j-c_\lambda\|^2 \leq a^2
(v_{nt_n^{-1}}(\tau_j))^2\right)&\leq& \pi_j\left(\inf_{\lambda \in
\Lambda/\Lambda^\prime} \|c_j-c_\lambda\|^2 \leq a^2
(v_{nt_n^{-1}}(\tau_j))^2\right)
\end{eqnarray*}
\begin{eqnarray}\label{proba_inf}
+\quad \pi_j\left(\inf_{\lambda \in \Lambda^\prime}
\|c_j-c_\lambda\|^2 \leq a^2 (v_{nt_n^{-1}}(\tau_j))^2\right)
\end{eqnarray}
Consider the particular subset $\Lambda'$ defined by
$$\Lambda'=\{\lambda \in  \Lambda: \; \| c_{\lambda_0} - c_\lambda\|^2 \leq  6 C_1^2  M^{2}
 (v_{nt_n^{-1}}(\tau_j))^2\}.$$
Notice that
$$
\lambda \in  \Lambda/\Lambda'\Longrightarrow \| c_\lambda -c_j\|^2
\geq a^2(v_{nt_n^{-1}}(\tau_j))^2
$$ due to the choice of $C_1$. It implies that the first term in the right hand side of (\ref{proba_inf}) is null and then, it
remains to prove that
\begin{eqnarray}
\lim_{n \rightarrow \infty} \pi_j\left(\inf_{\lambda \in \Lambda'}
\|c_j-c_\lambda\|^2 \leq a^2 (v_{nt_n^{-1}}(\tau_j))^2\right) =0.
\label{lambda'}
\end{eqnarray}
Since  $\lambda$ is in $\Lambda^\prime$, we get
 \begin{eqnarray*}
 \| c_\lambda -c_j\|_2^2  & =& \|
c_{\lambda_0} - c_\lambda\|^2 + \sum_k \sum_{\epsilon} u_{j}(n)^2 +
2 \sum_k\sum_{\epsilon} \delta_{k} u_{j}(n)
B_{j,k,\lambda,\lambda_0}
\\
& \geq & 3 C_1^2  M^{2}   \; (v_{nt_n^{-1}}(\tau_j))^2 + 2 \sum_k
\delta_{k} u_{j}(n)  \sum_{\epsilon} B_{j,k,\lambda,\lambda_0},
\end{eqnarray*}
where $$B_{j,k,\lambda,\lambda_0}= \int
\psi^{\epsilon}_{j,k}(c_{\lambda_0} - c_\lambda).$$
 Therefore assertion
(\ref{lambda'}) is equivalent to
\begin{eqnarray}
\lim_{n \rightarrow \infty} \pi_j\left(\inf_{\lambda \in \Lambda'}
2\sum_k \delta_{k}u_{j}(n) B_{j,k,\lambda,\lambda_0}\leq -a^2
(v_{nt_n^{-1}}(\tau_j))^2\right)=0 \nonumber.
\end{eqnarray}
or
\begin{eqnarray}
\lim_{n \rightarrow \infty} \pi_j\left(\sup_{\lambda \in \Lambda'}
2\sum_k (-\delta_{k})u_{j}(n) B_{j,k,\lambda,\lambda_0}\geq a^2
(v_{nt_n^{-1}}(\tau_j))^2\right)=0 \label{reseau}.
\end{eqnarray}

We can construct in the Euclidean metric an $\eta$-net $\Lambda'_{\eta}$  on the subset $\Lambda'$. For any $\lambda$ in $\Lambda'$, denote ${\tilde \lambda}$ the closest element in $\Lambda'_{\eta}$ to $\lambda$ in the Euclidean sense. Then for any  $\lambda \in \Lambda'$, we have by assumption {\bf A0}:
\begin{eqnarray*}
|\sum_k \delta_{k}u_{j}(n)( B_{j,k,\lambda,\lambda_0}- B_{j,k,{\tilde \lambda},\lambda_0})| & \leq & u_{j}(n)\sum_k |B_{j,k,\lambda,\lambda_0}- B_{j,k,{\tilde \lambda},\lambda_0}| \\
& \leq & u_{j}(n)\sum_k Q \eta^{\nu} 2^{-j} \| \psi^{\epsilon}\|_{\infty}\\
&\leq & 2^{j(s+1)} 2^{2j}  Q \eta^{\nu} 2^{-j} \| \psi^{\epsilon}\|_{\infty}\\
& \leq & \kappa 2^{j s}  \eta^{\nu},
\end{eqnarray*}  where $\kappa$ is a positive constant depending on $Q$, $C_1$, $M$ and $\| \psi^{\epsilon}\|_{\infty}$.
Chosing $\eta=n^{-b}$, with $b\nu  > \frac{s_{{\rm max}}}{ 2 s_{{\rm max}} + 1}$, the proof of relation  (\ref{reseau}) is then reduced to the proof of
\begin{eqnarray}
\lim_{n \rightarrow \infty} {\rm Card} (\Lambda'_{n^{-b}})  \pi_j\left(
2\sum_k (-\delta_{k})u_{j}(n) B_{j,k,{\tilde \lambda},\lambda_0}\geq a^2
(v_{nt_n^{-1}}(\tau_j))^2\right)=&0& \nonumber \\
\lim_{n \rightarrow \infty} (T n^{b})^{d_{\Lambda}}  \pi_j\left(
2\sum_k (-\delta_{k})u_{j}(n) B_{j,k,{\tilde \lambda},\lambda_0}\geq a^2
(v_{nt_n^{-1}}(\tau_j))^2\right)&=&0, \label{Bern-reseau}
\end{eqnarray}
where $\mbox{\rm {Diam}}$ is the diameter of $\Lambda$.
 \noindent
Finally, relation (\ref{lambda'}) is proved applying Bernstein
inequality in the right hand side of relation (\ref{Bern-reseau}).
Indeed Bernstein inequality is applied to
$$\displaystyle{\pi_j\left( 2\sum_k (-\delta_{k})u_{j}(n)
B_{j,k,{\tilde \lambda},\lambda_0}\geq a^2
(v_{nt_n^{-1}}(\tau_j))^2\right)},$$ with the   i.i.d. centered
random variables $Z_k= -\delta_{k} B_{j,k,{\tilde
\lambda},\lambda_0}$. In particular, Notice that $|Z_k| <
K_1v_{nt_n^{-1}}(\tau_j)$,  $\sum_k \V (Z_k) \leq K_2
(v_{nt_n^{-1}}(\tau_j))^2$, where $K_1$  and $K_2$ are positive
constants. Notice also that  it leads to an exponential bound of
order $\exp (- 2 ^{j})$.
\subsection{Proof of Relation (\ref{vrais})}
Set $$l_{n,\pi}= \displaystyle{\int \prod_{i=1}^n
\frac{c_j(U_i,V_i)}{c_{\lambda_0}(U_i,V_i)} d\pi_j(c_j)}.$$ Due to
the fact that the functions $\psi^{\epsilon}_{j,k}$ have disjoint
support, it is possible to rewrite $c_j$ as follows
$$
c_j=c_{\lambda_0} \prod_{k} (1 + \delta_{k} D_{j,k})
$$
for
$$
D_{j,k}=  u_{j}(n)
\sum_{\epsilon}\frac{\psi^{\epsilon}_{j,k}}{c_{\lambda_0}}.
$$
Then,
\begin{eqnarray*}
l_{n,\pi}& =& \prod_k \int \prod_{i=1}^n (1 + \delta_{k}
D_{j,k}(U_i,V_i)) d\pi_j (\delta_{k})\\& =& \prod_k \frac 1 2
\left\{\prod_{i=1}^n (1 +  D_{j,k}(U_i,V_i))+ \prod_{i=1}^n  (1 -
D_{j,k}(U_i,V_i))\right\},
\end{eqnarray*}
and
\begin{eqnarray*}
l_{n,\pi}^2 & =& \prod_k \frac 1 4  \left\{2 \prod_{i=1}^n \left[1 +
D^2_{j,k}(U_i,V_i)\right]
 + 2 \prod_{i=1}^n \left[1 - D^2_{j,k}(U_i,V_i)\right]\right. \\
 &&\left.+ H\left(D_{j,k}(U_i,V_i),\left(D^{b_t}_{j,k}(U_t,V_t)
 \right)_{t \in \{1,\ldots,i-1,i+1,\ldots,n\}}\right)\right\},
\end{eqnarray*}
where $b_t$ is either $0$ or $2$. Due to the independence of the
data and acting as in \citep{Pouet:2000}, it can be shown that
$$\E_{\lambda_0}\left [H\left(D_{j,k}(U_i,V_i),\left(D^{b_t}_{j,k}(U_t,V_t)
 \right)_{t \in \{1,\ldots,i-1,i+1,\ldots,n\}}\right)\right]=0.$$
Therefore,
\begin{eqnarray*}
\E_{\lambda_0} [l_{n,\pi}^2(U_i,V_i)] & \leq & \prod_k \left\{
\left(1 + \E_{\lambda_0}D_{j,k}^2(U_i,V_i)\right)^n
+\left(1 - \E_{\lambda_0}D_{j,k}^2(U_i,V_i)\right)^n \right\}\\
& \leq&  \prod_k \cosh \left( n
\E_{\lambda_0}D_{j,k}^2(U_i,V_i)\right).
\end{eqnarray*}
Using the inequality $\log (\cosh(u)) \leq K u^2$ where $K$ is a
fixed constant and since $c_{\lambda_0}$ is bounded from below by
$m$, one obtains
\begin{eqnarray*}
\frac{1}{N_n^2} \sum_{j \in {\tilde J}} \exp(\log (\E_{\lambda_0}
l_{n,\pi})^2  ) &\leq & \frac{1}{N_n^2} \sum_{j \in {\tilde J}} \exp \left\{  K
n^2
 \sum_k \left( \E_{\lambda_0}D_{j,k}^2(U_i,V_i)\right)^2\right\} \\
&\leq & \frac{1}{N_n^2} \sum_{j \in {\tilde J}} \exp \left\{  \frac{3^2 K}{m^2}n^2 2^{2j} u_{j}(n)^4\right\}\\
&\leq& \frac{\log(n)^{\kappa}}{\log(n)(1+o_n(1)) },
\end{eqnarray*}
where $\kappa =  K (3 C_1^2 M^2)^{2} m^{-2}=4Ka^4m^{-2}$. Choosing
$a$ small enough and $\kappa<1$, Relation (\ref{vrais}) is then
proved.

\section{Appendix A: Proof of Lemma \ref{etude-stat}}
In this part, $\kappa$ denotes any positive constant which may
depend on $\phi$, $M$ and on $\|c\|,\|c_\lambda\|$.

\subsection{Notations and Preliminaries} Let us  define or recall some notations that
will be used below. For any $k \in \BBz^2$, set
\begin{eqnarray*}
\xi_k(X_i,Y_i) & =&
\phi_{j,k}\left(\widehat{F}(X_i),\widehat{G}(Y_i) \right)-
\phi_{j,k}\left(F(X_{i}),G(Y_{i})\right) \label{xi} \\
\omega_{j,k}^\lambda(X_{i},Y_{i}) & =& \phi_{j,k}(F(X_{i}),
G(Y_{i}))-c_{j,k}(\lambda) \nonumber \\
 \omega_{j,k}^\infty(X_{i},Y_{i}) & =& \phi_{j,k}(F(X_{i}),
G(Y_{i}))-c_{j,k}, \nonumber
\end{eqnarray*}
where $i$ is in ${\cal I}_2$. First,  the localization property of the
scaling function implies that only few $\xi_k(X_i,Y_i)$ will be used since the others are zero.
Indeed, one has  the following result

\vspace{0.2cm}

\begin{lemma}\label{nombredei}For any $k \in \BBz^2$, let us denote
\begin{eqnarray*}
N_{j}&=& \mbox{ {\rm card } }\left\{i \in {\cal I}_2; \xi_k(X_i,Y_i)\not=
0\right\}\label{ni}.
\end{eqnarray*}
Let $\delta>0$. For any level $j$ such that
$$2^j\leq
\frac{2}{3\sqrt{\delta+1}}\,\left(\frac{n_2}{\log(n_2)}\right)^{1/2},$$
 one has
\begin{eqnarray*}\label{nx}
\mathbb{P}(N_{j}>2(2L+3)n_22^{-j}))&\leq
&K(n_1^{-\delta}+n_2^{-\delta}).
\end{eqnarray*}
\end{lemma}
We refer to  \citep{Genest/Masiello/Tribouley:2008} for the proof of
this lemma since a similar result is established with an estimate
$\widehat{F}$ built on the whole sample: it  guarantees
 in particular that
 $\widehat{F}(X_{(i:n)})=i/n$, where $X_{(i:n)}$ denotes the $i-$th (among $n$) order statistic.
 In our case, the situation is different since $\widehat{F}(X_{(i:n_2)})$ is based
 on the observations lying in the subsample whose indices are in ${\cal I}_1$ whereas  it is calculated in an observation
  lying in the subsample whose indices are in ${\cal I}_2$ ; nevertheless, applying the Dvoretsky--Kiefer--Wolfovitz Inequality,
 the following deviation inequality  holds. For any $\epsilon>0$,
$\mathbb{P}_{\widehat{F}}=\mathbb{P}\left(\left|\widehat{F}(X_{(i:n_2)})-\frac{i}{n_2}\right|\geq
2\epsilon\right)$ is bounded from above by

\vspace{0.3cm}
\begin{eqnarray*}
\mathbb{P}_{\widehat{F}}&\leq&
\mathbb{P}\left(\left|\widehat{F}(X_{(i:n_2)})-{F}(X_{(i:n_2)})\right|\geq
\epsilon\right)+
\mathbb{P}\left(\left|{F}(X_{(i:n_2)})-\widetilde{F}(X_{(i:n_2)})\right|\geq
\epsilon\right)\\
&\leq& \mathbb{P}\left(\|\widehat{F}-{F}\|_\infty\geq
\epsilon\right)+ \mathbb{P}\left(\|\widetilde{F}-{F}\|_\infty\geq
\epsilon\right)\\
&\leq& K\left(n_1^{-\delta}+n_2^{-\delta}\right),
\end{eqnarray*}
as soon as we take $\epsilon=\sqrt{\delta\log(n_1)/(2n_1)}\vee
\sqrt{\delta\log(n_2)/(2n_2)}$. Here $\widehat{F}$ represents the
empirical margin computed with the  subsample whose indices in ${\cal I}_1$  and
$\widetilde{F}$, the empirical margin computed with the  subsample whose indices in ${\cal I}_2$.

\subsubsection{Study of  $\widehat{T_j}(\lambda)$}
Rewrite $\widehat{\theta_{j,k}}(\lambda)$ in
$\widehat{T_j}(\lambda)=\sum_k \widehat{\theta_{j,k}}(\lambda)$ as
follows
\begin{eqnarray*}
\widehat{\theta_{j,k}}(\lambda)&=&
\frac{1}{n_2(n_2-1)}\mathop{\sum_{i_1, i_2 \in {\cal I}_2}}_{ i_1 \neq i_2
}\omega_{j,k}^\lambda(X_{i_1},Y_{i_1}) \omega_{j,k}^\lambda(X_{i_2},Y_{i_2}).
\end{eqnarray*}
For all $i \in {\cal I}_2$, one has $\E
(\omega_{j,k}^\lambda(X_{i},Y_{i}))= c_{j,k} -c_{j,k}(\lambda)$, which implies that
\begin{eqnarray*}
\E ( \widehat{T_j}(\lambda)) = \sum_k \theta_{j,k}(\lambda)=
T_j(\lambda).\end{eqnarray*} Moreover for $p \neq k$, one obtains

 \noindent
$\E
(\widehat{\theta_{j,k}}(\lambda)\widehat{\theta_{j,p}}(\lambda))$
{\footnotesize
\begin{eqnarray*}
&= & \frac{1}{(n_2(n_2-1))^2 }\sum_{i_1\not =i_2\not=i_3\not=i_4} \E
\left[\omega_{j,k}^\lambda(X_{i_1},Y_{i_1})\right] \E
\left[\omega_{j,p}^\lambda(X_{i_3},Y_{i_3})\right]
\E \left[\omega_{j,k}^\lambda(X_{i_2},Y_{i_2})\right]\E \left[\omega_{j,p}^\lambda(X_{i_4},Y_{i_4})\right]\\
& & +4\frac{1}{(n_2(n_2-1))^2  }\sum_{i_1\not =i_2\not=i_3} \E
\left[\omega_{j,k}^\lambda(X_{i_1},Y_{i_1})\right] \E
\left[\omega_{j,p}^\lambda(X_{i_3},Y_{i_3})\right]
\E \left[\omega_{j,k}^\lambda(X_{i_2},Y_{i_2})\omega_{j,p}^\lambda(X_{i_2},Y_{i_2})\right]\\
& & +2\frac{1}{(n_2(n_2-1))^2 }\sum_{i_1\not =i_2} \E
\left[\omega_{j,k}^\lambda(X_{i_1},Y_{i_1})\omega_{j,p}^\lambda(X_{i_1},Y_{i_1})\right]
\E
\left[\omega_{j,k}^\lambda(X_{i_2},Y_{i_2})\omega_{j,p}^\lambda(X_{i_2},Y_{i_2})\right]\\
&\leq &\theta_{j,k}(\lambda)\theta_{j,p}(\lambda)+
\frac{4}{n_2}\,\left(c_{j,k} -c_{j,k}(\lambda)\right)\left(c_{j,p}
-c_{j,p}(\lambda)\right) \left[\int\left(\phi_{j,k}-\int
\phi_{j,k}c_\lambda\right)\left(\phi_{j,p}-\int
\phi_{j,p}c_\lambda\right)c\right]\\&&\,+
\frac{2}{n_2(n_2-1)}\left[\int\left(\phi_{j,k}-\int
\phi_{j,k}c_\lambda\right)\left(\phi_{j,p}-\int
\phi_{j,p}c_\lambda\right)c\right]^2,\end{eqnarray*}} which implies
that
{\footnotesize
\begin{eqnarray*}
\V (\widehat{T_j}(\lambda))&=
&\E\left(\sum_{k}\widehat{\theta_{j,k}}(\lambda)\right)^2-\left(\E\sum_{k}\widehat{\theta_{j,k}}(\lambda)\right)^2
\\
&\leq& \frac{4}{n_2}\,\sum_{k,p}\left(c_{j,k}
-c_{j,k}(\lambda)\right)\left(c_{j,p} -c_{j,p}(\lambda)\right)
\left[\int\left(\phi_{j,k}-\int
\phi_{j,k}c_\lambda\right)\left(\phi_{j,p}-\int
\phi_{j,p}c_\lambda\right)c\right]\\&&\,+
\frac{2}{n_2(n_2-1)}\sum_{k,p}\left[\int\left(\phi_{j,k}-\int
\phi_{j,k}c_\lambda\right)\left(\phi_{j,p}-\int
\phi_{j,p}c_\lambda\right)c\right]^2.
 \end{eqnarray*}}
  Applying the
H\"older Inequality and the consequence of the Parseval Equality, we
get

\vspace{0.2cm}

{\footnotesize
 $\displaystyle{
 \sum_{k
  p}\left[\int\left(\phi_{j,k}-\int
\phi_{j,k}c_\lambda\right)\left(\phi_{j,p}-\int
\phi_{j,p}c_\lambda\right)c\right]^2}$
\begin{eqnarray*}& \leq & 2^2 \left(
\sum_{k,p}\left[\int\phi_{j,k}\phi_{j,p}c\right]^2+
2\sum_{k,p}\left[\int\phi_{j,k}c_\lambda\int\phi_{j,p}c\right]^2
+\left(\sum_{k}\left[\int\phi_{j,k}c_\lambda\right]^2\right)^2
\right)
\\
&\leq & 2^2  \left(  \left( \sum_{k} \int\phi_{j,k}^2c \right)^2 +2
\int c^2   \int c_\lambda^2 +\left(  \int c_\lambda^2\right)^2
 \right)\leq  \kappa \,2^{2j}.
\end{eqnarray*}}
We conclude that
\begin{eqnarray*}
\V (\widehat{T_j}(\lambda))&\leq & \kappa \left(
\frac{4}{n_2}\left(\sum_{k,p}\theta_{j,k}(\lambda)\theta_{j,p}
(\lambda)\right)^{1/2}\;2^{j}+\frac{2^{2j}}{n_2(n_2-1)}\right)\\
&\leq & \kappa \left( \frac{2^j}{n_2}T_j(\lambda)+
\,\frac{2^{2j}}{n_2^2}
 \right),
\end{eqnarray*}
which is the announced result for $\widehat{T_j}(\lambda)$.

\subsubsection{Study of  $T_j^\spadesuit $ and $T_j^{\clubsuit }(\lambda)$    }
 Let us denote
\begin{eqnarray*}
 A_{i_1}&=&
\left[\xi_{k}(X_{i_1},Y_{i_1})\right],\quad D_{i_1}=
\sum_{k,p}\left(\E
\left[\xi_{k}(X_{i_1},Y_{i_1}) \xi_{p}(X_{i_1},Y_{i_1})\right]\right)^2\\
B_{i_1,i_2}&=&\sum_{k}\left[\xi_{k}(X_{i_1},Y_{i_1})
\xi_{k}(X_{i_2},Y_{i_2}\right],\; C_{i_1,i_2}=
\xi_{k}(X_{i_1},Y_{i_1}) \xi_{p}(X_{i_2},Y_{i_2}).
\end{eqnarray*}
We need the following results which are stated in the lemma below
\begin{lemma}\label{horreur}
Assume that the scaling function is $q$-differentiable. For any
level $j\leq j_\infty$, there exists some positive constant $\kappa$
depending on $\phi$, its derivatives and on $\|c\|_\infty$ (which
might be $\|c_\lambda\|_\infty$ for some $\lambda \in \Lambda$) such
that for any distinct indices $i_1,i_2$, one obtains
\begin{eqnarray}\label{A}
\E |A_{i_1}|&\leq &\kappa\left(\frac{\log(n_1)}{n_1}\right)^{1/2}
\end{eqnarray}
\begin{eqnarray}\label{BCD}
\E |B_{i_1,i_2}| & \leq & \kappa
2^{2j}\left(\frac{\log(n_1)}{n_1}\right),\quad \E |C_{i_1,i_2}|\leq
\kappa \left(\frac{\log(n_1)}{n_1}\right)\\
|D_{i_1}|& \leq &
2^{6j}\left(\frac{\log(n_1)}{n_1}\right)^2.\nonumber
\end{eqnarray}
\end{lemma}
We prove relation (\ref{A}) in the next section, relations
(\ref{BCD}) are  proven
 in  \citep{Genest/Masiello/Tribouley:2008}. We have
\begin{eqnarray*}
\E T_{j}^\spadesuit  &=&\frac{1}{n_2(n_2-1)}\mathop{\sum_{i_1, i_2
\in {\cal I}_2 }}_{ i_1 \neq i_2 }\E[B_{i_1,i_2}].
\end{eqnarray*}
Using Lemma \ref{nombredei} and Lemma \ref{horreur}, it follows
\begin{eqnarray*}
\E |T_{j}^\spadesuit | &\leq &\frac{1}{n_2(n_2-1)}\,(n_22^{-j})^2\,
2^{2j}\left(\frac{\log(n_1)}{n_1}\right)\leq\left(\frac{\log(n_1)}{n_1}\right). \end{eqnarray*}
Moreover, we get
\begin{eqnarray*}
 T_{j}^\clubsuit (\lambda) &=&\frac{1}{n_2(n_2-1)}\mathop{\sum_{i_1,
i_2 \in {\cal I}_2}}_{ i_1 \neq i_2 }\sum_k
\left[\xi_k(X_{i_1},Y_{i_1})\omega_{j,k}^\lambda(X_{i_2},Y_{i_2})\right].
\end{eqnarray*}
By H\"older Inequality and from  lemmas \ref{nombredei} and
\ref{horreur}, one obtains
\begin{eqnarray*}
\E| T_{j}^\clubsuit (\lambda)|&\leq
&\frac{1}{n_2(n_2-1)}\mathop{\sum^{n_2}_{i_1, i_2 \in {\cal I}_2}}_{
i_1 \neq i_2 }\left(\sum_k\,(\E(A_i))^2\;\sum_k(
\E\omega_{j,k}^\lambda(X_{i_2},Y_{i_2}))^2\right)^{1/2}.\end{eqnarray*}
Remembering that $\E
\omega_{j,k}^\lambda(X_{i_2},Y_{i_2})=(c_{j,k}-c_{j,k}(\lambda))$
for any index $i_2$, we get
\begin{eqnarray*}
\E| T_{j}^\clubsuit (\lambda)|&\leq&
\frac{1}{n_2(n_2-1)}\,(n_22^{-j})n_2\,
\left[2^{2j}\frac{\log(n_1)}{n_1}\,T_j(\lambda)\right]^{1/2}\\
&\leq&K\,\left(\frac{\log(n_1)}{n_1}\,T_j(\lambda)\right)^{1/2}.
\end{eqnarray*}
Let us study the moments of $T_{j}^\clubsuit   (\lambda)$ under
$\P_\lambda$. Since $\E_\lambda\omega_{j,k}^\lambda(X_i,Y_i)=0$ for
any $k$ and $i$, we obviously have $\E_\lambda T_{j}^\clubsuit
(\lambda)=0$ and
$$\E_\lambda (T_{j}^\clubsuit (\lambda))^2=\left(\frac{1}{n_2(n_2-1)}\right)^2\sum_{i_1\not=i_2
}T_{i_1,i_2}+\left(\frac{1}{n_2(n_2-1)}\right)^2\sum_{i_1\neq i_2\neq
i_3 }S_{i_1,i_2,i3},$$ where
\begin{eqnarray*}
T_{i_1,i_2}&=&\sum_{k,p}\left(\E_\lambda\left[\xi_k(X_{i_1},Y_{i_1})\xi_p(X_{i_1},Y_{i_1})\right]
\E_\lambda \left[\omega_{j,k}^\lambda(X_{i_2},Y_{i_2})\omega_{j,p}^\lambda(X_{i_2},Y_{i_2})\right]\right),\\
S_{i_1,i_2,i_3}&=& \sum_{k,p}\left(\E_\lambda
\left[\xi_{k}(X_{i_1},Y_{i_1})
\xi_{p}(X_{i_2},Y_{i_2})\right]\E_\lambda\left[\omega_{j,k}^\lambda(X_{i_3},Y_{i_3})\omega_{j,p}^\lambda(X_{i_3},Y_{i_3})\right]\right).
\end{eqnarray*}
By H\"older Inequality, we have
\begin{eqnarray*}
T_{i_1,i_2}&=&\sum_{k,p}\E_\lambda
\left[\xi_k(X_{i_1},Y_{i_1})\xi_p(X_{i_1},Y_{i_1})\right] \E_\lambda
\left[\omega_{j,k}^\lambda(X_{i_2},Y_{i_2})\omega_{j,p}^\lambda(X_{i_2},Y_{i_2})\right]\\&\leq&D_{i_1}^{1/2}
\left(\sum_{k,p}\left( \E_\lambda
\left[\omega_{j,k}^\lambda(X_{i_2},Y_{i_2})\omega_{j,p}^\lambda(X_{i_2},Y_{i_2})\right]\right)^2\right)^{1/2}
\end{eqnarray*}
With Parseval Equality, we get
\begin{eqnarray*}
\sum_{k,p}\left(\E_\lambda
\left[\omega_{j,k}^\lambda(X_{i_2},Y_{i_2})\omega_{j,p}^\lambda(X_{i_2},Y_{i_2})\right]\right)^2&\leq&
\sum_{k,p}\left(\int\phi_{j,k}\phi_{j,p}c_\lambda\right)^2\\
&\leq& K\,\sum_{k}\int\phi^2_{j,k}c_\lambda^2\leq K\,2^{2j},
\end{eqnarray*}
which  combining with Lemma \ref{horreur}, implies that
\begin{eqnarray*}
T_{i_1,i_2}&\leq&K\,\left(2^{6j}\left(\frac{\log(n_1)}{n_1}\right)^2\right)^{1/2}
\left(2^{2j}\right)^{1/2}\leq
2^{4j}\left(\frac{\log(n_1)}{n_1}\right).
\end{eqnarray*}
In the same way,
\begin{eqnarray*}
S_{i_1,i_2,i_3}&\leq &K\,\left(\sum_{k,p}(\E_\lambda
C_{i_1,i_2})^2\right)^{1/2} \left(\sum_{k,p}\left(
\E_\lambda\left[\omega_{j,k}^\lambda(X_{i_2},Y_{i_2})\omega_{j,p}^\lambda(X_{i_2},Y_{i_2})\right]\right)^2\right)^{1/2}
\\&\leq&
\left(2^{2j}\right)^{1/2}\left(2^{4j}\left(\frac{\log(n_1)}{n_1}\right)^2\right)^{1/2}
\leq 2^{3j}\left(\frac{\log(n_1)}{n_1}\right).
\end{eqnarray*}
From  Lemma \ref{nombredei}, one has
\begin{eqnarray*}
\E_\lambda( T_{j}^\clubsuit  (\lambda))^2&\leq&
K\;\frac{1}{n_2^2(n_2-1)^2}(n_22^{-j})n_22^{4j}\left(\frac{\log(n_1)}{n_1}\right)\\
&&+
K\;\frac{1}{n_2^2(n_2-1)^2}(n_22^{-j})^2n_22^{3j}\left(\frac{\log(n_1)}{n_1}\right)\\
&\leq&K\,2^{j}\left(\frac{\log(n_1)}{n_2n_1}\right).
\end{eqnarray*}

\subsection{ Proof of Lemma \ref{horreur}}
The following expansion is crucial because it allows to reduce the
study to univariate variables.
\begin{eqnarray}\label{biuni}
&& \xi_{k}(X_i,Y_i)= \xi_{k_1}(X_i)\xi_{k_2}(Y_i)\\&&\hspace{2cm}+
 \xi_{k_1}(X_i)\phi_{jk_2}\left(G(Y_i)\right)
+
 \xi_{k_2}(Y_i)\phi_{jk_1}\left(F(X_i)\right)
,\nonumber \end{eqnarray} where the univariate statistics
$\xi_{k_1}(X_i)$ and $\xi_{k_2}(Y_i)$ are defined as follows
\begin{eqnarray*}\label{xiuni}\xi_{k_1}(X_i)&=&\phi_{j,k_1}\left(\frac{\widehat{F}(X_i)}{n_1}\right)-\phi_{j,k_1}(F(X_i))\\
\xi_{k_2}(Y_i)&=&\phi_{j,k_2}\left(\frac{\widehat{G}(Y_i)}{n_1}\right)-\phi_{j,k_1}(G(Y_i)).
\end{eqnarray*}
 Assuming that $\phi$ is continuously
$q-$differentiable, we get
\begin{eqnarray*}\label{taylor}\xi_{k_1}(X_i)&=&\hat z_{k_1}(X_i)+\hat
w_{k_1}(X_i),
\end{eqnarray*} where
\begin{eqnarray*}
 \hat
z_{k_1}(X_i)&=&\sum_{\ell=1}^{q-1}\frac{2^{j\ell}}{\ell!}(\widehat{F}(X_i)-F(X_i))^\ell\;\phi^{(\ell)}_{j,k_1}(F(X_i))
\end{eqnarray*} and
\begin{eqnarray*}
\hat w_{k_1}(X_i)&=&2^{qj}\int_{\widehat{ F}(X_i)}^{F(X_i)}
\phi^{(q)}_{j,k_1}(t)\;( F(X_i)-t)^{q-1}dt.
\end{eqnarray*}
A direct application of the Dvoretsky, Kiefer and
 Wolfovitz Inequality leads to the following bound
\begin{eqnarray*}\label{DKW}
\P(\|\widehat{ F}-F\|_\infty
>\epsilon)\leq K\,\exp(-2n_1\epsilon^2)\leq Kn_1^{-\delta},
\end{eqnarray*}
as soon as $\epsilon= \sqrt{0.5\;\delta\log(n_1)/n_1}$. In the
sequel, we take such an $\epsilon$ with $\delta$ large enough. Since
$j\leq j_\infty$ where $j_\infty$ is defined in (\ref{j0jinfty}),
observe that $2^j\epsilon\leq 1$ and then we get
\begin{eqnarray*}
|\hat z_{k_1}(X_i)|&\leq &K\,2^{j}\epsilon\max_{\ell=1,\ldots
q-1}\,|\phi^{(\ell)}_{j,k_1}(F(X_i))|(1+o_P(1))\\
|\hat w_{k_1}(X_i)|&\leq & K\,2^{(q+1/2)j}\epsilon^{q}(1+o_P(1))
\end{eqnarray*}
which leads to the following bound
{\footnotesize
\begin{eqnarray*}\label{xibound1}|\xi_{k_1}(X_i)|&\leq &K\left(
\;2^{(q+1/2)j}\epsilon^q+2^j\epsilon\max_{\ell=1,\ldots q-1}
\,|\phi^{(\ell)}_{j,k_1}(F(X_i))|\right)(1+o_P(1)).
\end{eqnarray*}}
 The same kind of
result obviously holds for $\xi_{k_2}(Y_i)$. In the sequel, we
need the following evaluations (which also hold for any derivatives
of $\phi$).
 Using expansion (\ref{biuni}), we get
\begin{eqnarray*}
 \xi_k(X_i,Y_i)&=&S_1+S_2,
\end{eqnarray*}
where
\begin{eqnarray*}
S_1&=& \xi_{k_1}(X_i)\xi_{k_2}(Y_i), \\S_2&=&
\xi_{k_1}(X_i)\phi_{j,k_2}\left(G(Y_i)\right)+
\xi_{k_2}(Y_i)\phi_{j,k_1}\left(F(X_i)\right).
 \end{eqnarray*}
Using (\ref{xibound1}), we get
\begin{eqnarray*}
\E |S_1|&\leq &K\,\left( 2^{(2q+1)j}\epsilon^{2q}+2^{(q+1)j}\epsilon^{q+1}+2^j\epsilon^2\right),\\
\E |S_2|&\leq&K\,\left(2^{qj}\epsilon^{q}+\epsilon\right).
 \end{eqnarray*}
If  $2^j\leq (n_1/\log(n_1))^{1/2-1/2q}$, we obtain $
\E|\xi_k(X_i,Y_i)|\leq\epsilon$ which ends the proof.

\section{Appendix B : Proof of Lemma \ref{Bernstein}}
Let us denote
 $T_{j}^\diamond(\lambda)=n_2^{-1}\sum_{i \in {\cal I}_2} Z_i$ where
\begin{eqnarray*}
Z_i&=&\sum_k\left( \phi_{jk}(F(X_i),G(Y_i))-c_{jk}\right)\left(
c_{jk}-c_{jk}(\lambda)\right),\\
\E_cZ_i&=&0,\\
 |Z_i|&\leq&\left(\sum_k\left(
\phi_{jk}(F(X_i),G(Y_i))-c_{jk}\right)^2 \sum_k \left(
c_{jk}-c_{jk}(\lambda)\right)^2\right)^{1/2}\\
&\leq&K_1\,\left(\left( 2^j\right)^2 T_j(\lambda)\right)^{1/2}\leq
K\,2^{j}
T_j(\lambda)^{1/2},
\end{eqnarray*}
and
\begin{eqnarray*} \V_c(Z_i)&\leq& \sum_{k,p} \E\left(
\phi_{jk}(F(X_i),G(Y_i))-c_{jk}\right)\left(\phi_{jp}(F(X_i),G(Y_i))-c_{jp}\right)\\
&&\hspace{1 cm}\times \left|\left(
c_{jk}-c_{jk}(\lambda)\right)\left(
c_{jp}-c_{jp}(\lambda)\right)\right|\\
&\leq& \sum_{k,p}\left( \E
\phi_{jk}^2(F(X_i),G(Y_i)) \E\phi_{jp}^2(F(X_i),G(Y_i))\right)^{1/2}\\
&&\hspace{1 cm}\times \left|\left(
c_{jk}-c_{jk}(\lambda)\right)\left(
c_{jp}-c_{jp}(\lambda)\right)\right|\\
&\leq& \|c\|_\infty\,\;\sum_{k}\left(
c_{jk}-c_{jk}(\lambda)\right)^2=\|c\|_\infty\, T_j(\lambda).
\end{eqnarray*}
 Applying
Bernstein Inequality to the $Z_i'$ s leads to  prove Lemma
\ref{Bernstein}.

\section{Appendix C: Proof of Lemma \ref{LDUstat}}
\subsection{$U$-Statistic}
Let us first recall  the result of \citep{Gine/Latala/Zinn:2000}.
\begin{proposition} {\sc (Theorem 3.3 p. 21 \citep{Gine/Latala/Zinn:2000})}\\
It exists an universal positive constant $K_g < \infty$ such that,
if $\Omega$ is a bounded canonical kernel of two variables for the
i.i.d. $Z_{i_1},Z_{i_2}$, $i_1, \; i_2 \in \{1,\ldots,{\tilde n}\}$,
where ${\tilde n} \in \BBn $, for any $x>0$, we have
\begin{eqnarray*}
\P \left( |\sum_{i_1,i_2} \Omega(Z_{i_1},Z_{i_2})| > x \right) \leq
K_g \exp \left( - \frac 1 K_g \; \min \left\{\frac{x^2}{C^2},
\frac{x}{D},\left(\frac{x}{B}\right)^{2/3},
\left(\frac{x}{A}\right)^{1/2} \right\}\right), \label{Gine}
\end{eqnarray*}
where
\begin{eqnarray*}
A &= &\|\Omega (\cdot,\cdot)\|_{\infty},\, B^2= \tilde n \left[  \|
\E[ \Omega^2(Z_1,\cdot)] \|_{\infty}+ \|  \E [\Omega^2(\cdot,Z_2)]
\|_{\infty} \right],\\
 C^2& =& \tilde n^2 \E [
(\Omega(Z_1,Z_2))^2] \end{eqnarray*} and
\begin{eqnarray*}
D= \tilde n \sup_{\Omega_1,\Omega_2} \{ \E[\Omega(Z_{1},Z_{2})
\Omega_1(Z_{1}) \Omega_2(Z_{2})]:
 \E [\Omega_1^2(Z_{1})]\leq 1;  \E [\Omega_2^2(Z_2)] \leq 1 \}.
\end{eqnarray*}
\end{proposition}

\vspace{0.3cm}

 We apply this proposition for $Z_i=(F(X_i),G(Y_i))$,
$\tilde n=n_2$  and
 the kernel
\begin{eqnarray*}
\Omega_{{\tilde c}}\left( Z_{i_1}, Z_{i_2}\right) & =& \sum_k \left\{
\phi_{j,k}\left(Z_{i_1} \right) - \E_{{\tilde c}}[ \phi_{j,k}\left(Z_{i_1}
\right)]\right\} \times   \left\{ \phi_{j,k}\left(Z_{i_2} \right) -
\E_{{\tilde c}}[ \phi_{j,k}\left(Z_{i_2} \right)]\right\},
\end{eqnarray*}
which is considered  under the distribution $\P_{{\tilde c}}$ where ${\tilde c}$ is either $c_\lambda$ or $c$.
The quantities $A$, $B$, $C$ and $D$ are evaluated in the following
lemma which is proved in the next section.
\begin{lemma} \label{ABCDBornes}
There exists some positive constant $K_1$ larger than either
$$
\left(12 L^2\|\phi\|_{\infty}^2\right) \vee \left(
2\|{\tilde c}\|_\infty\right) \vee \left( 2 L^2
\|\phi\|_{\infty}^2\right)\vee \left(4
\|{\tilde c}\|_\infty(\|{\tilde c}\|_\infty+3L^4\|\phi\|_\infty^2))\right) $$
 such that
 $$A \leq K_1
2^{2j},\;B^2 \leq K_1 \,n_2 2^{2j},\;C^2 \leq K_1\,n_2^2 2^{2j},\;D \leq
K_1\, n_2  ,$$
where ${\tilde c}$ is either $c_\lambda$ or $c$.
\end{lemma}
\vspace{0.3cm}

 Again define ${\tilde c}$ as  $c_\lambda$ or $c$, then applying both the result of \citep{Gine/Latala/Zinn:2000} and
Lemma \ref{ABCDBornes}, for any level $j$ and any $x\geq
2^j((n_2-1)n_2)^{-1/2}\sqrt{\log(\log(n_2))}$, it immediately
follows that {\footnotesize
\begin{eqnarray*}
\P_{{\tilde c}}\left(|\frac{1}{n_2(n_2-1)}\mathop{\sum_{i_1, i_2 \in {\cal
I}_2}}_{ i_1 \neq i_2 } \Omega_{{\tilde c}}\left( Z_{i_1}, Z_{i_2}\right)|>
\mu x \right)&\leq& K_g\exp\left(-\delta\,\log(\log(n_2))\right).
\end{eqnarray*}}
which ends the proof of Lemma \ref{LDUstat}.

\subsection{Proof of Lemma \ref{ABCDBornes}}
Let us denote $(U,V)=(F(X),G(Y))$ any pair
of random variables
whose marginal distribution are both uniform on $[0,1]$. Denote
 ${\tilde c}$ the copula density which is  $c_\lambda$ or $c$; in the same spirit, the coefficients
 ${\tilde c}_{j,k}$ stand for $c_{j,k}(\lambda)$ or $c_{j,k}$.
Recall that
\begin{eqnarray*}
c_{j,k}(\lambda)& =&  \E_{\lambda} \left[ \phi_{j,k}(F(X),G(Y))
\right] =\int c_\lambda (u,v )\phi_{j,k}(u,v)du dv \nonumber.\\
c_{j,k}& =&  \E \left[ \phi_{j,k}(F(X),G(Y)) \right] = \int c (u,v
)\phi_{j,k}(u,v)du dv \nonumber.
\end{eqnarray*}
Notice that
\begin{eqnarray*}
\sum_{k,p} \E_{{\tilde c}} [
\phi_{j,k}(U_{i_1},V_{i_1})\phi_{j,p}(U_{i_1},V_{i_1})
]& \leq & 2^{2j}, \\
\sum_k (\E_{{\tilde c}}[ \phi_{j,k}(U,V )])^2 & =& \sum_k {\tilde c}_{j,k}^2 \leq \|{\tilde c}\|^2
\leq M.
\end{eqnarray*}
We get
\begin{eqnarray*}
 A & =& \| \sum_{k}\left(\phi_{j,k}(u_1,v_1) - \E_{{\tilde c}}[
\phi_{j,k}(U,V )] \right) \left(\phi_{j,k}(u_2,v_2) - \E_{{\tilde c}}[
\phi_{j,k}(U,V )] \right) \|_{\infty} \\
&\leq &  \| \sum_{k} \phi_{j,k}(u_1,v_1)\phi_{j,k}
(u_2,v_2)\|_{\infty} + 2 \|
 \sum_{k} \phi_{j,k}(u_1,v_1) \E_{{\tilde c}}[ \phi_{j,k}(U,V )] \|_{\infty}\\&&+
  \|
 \sum_{k} \left(\E_{{\tilde c}}[ \phi_{j,k}(U,V )] \right)^2\|_\infty
 \\
& \leq & L^2\,2^{2j} \|\phi\|_{\infty}^2 +
2L^2\|\phi\|_{\infty}\|{\tilde c}\|_2  2^{j}   + \|{\tilde c}\|_2^2 \leq K 2^{2j},
\end{eqnarray*}
where $K\geq 2 L^2 \|\phi\|_{\infty}^2 $ and
\begin{eqnarray*}
B^2& =& 2 n_2 \left \| \sum_{k,p}  \E_{\tilde c} \left[
\left(\phi_{j,k}(U_{i_1},V_{i_1}) - \E_{\tilde c}[ \phi_{j,k}(U,V )] \right)
\left(\phi_{j,p}(U_{i_1},V_{i_1}) - \E_{\tilde c}[ \phi_{j,p}(U,V) ] \right)
\right]\right.  \\
& & \times \left.\left(\phi_{j,k}(u_2,v_2)  - \E_{\tilde c}[ \phi_{j,k}(U,V )]
\right) \left(\phi_{j,p}(u_2,v_2) - \E_{\tilde c}[ \phi_{j,p}(U,V) ]
\right)\right\|_{\infty} \\
& \leq & 2n_2  \left  \| \sum_{k,p}  \left|\int \phi_{j,k} \phi_{j,p}
{\tilde c} -\int \phi_{j,k}\,{\tilde c}\int \phi_{j,p}\,{\tilde c}  \right|
\left(\phi_{j,k}(u_2,v_2)  - \E_{\tilde c}[
\phi_{j,k}(U,V)] \right) \right.\\
& & \times \left.  \left(\phi_{j,p}(u_2,v_2)  - \E_{\tilde c}[ \phi_{j,p}(U,V)]
\right)\right\|_{\infty}   \\
& \leq &   2n_2 (2\|{\tilde c}\|_\infty) \left[\left   \| \sum_{k,p}
\phi_{j,k}(u_2,v_2)\phi_{j,p}(u_2,z_2)\right\|_{\infty} + 2\left \|
\sum_{k,p} \phi_{j,k}(u_2,v_2)\E_{\tilde c}[ \phi_{j,k}(U,V
)]\right\|_{\infty}\right.\\&&\left. + \left \| \sum_{k,p}\E_{\tilde c}[
\phi_{j,k}(U,V )]\E_{\tilde c}[ \phi_{j,p}(U,V )]\right\|_{\infty}
\right]  \\
& \leq &  (4 n_2\|{\tilde c}\|_\infty) \left(  2^{2j}L^4 2
\|\phi\|_\infty^2+2L^22^{j}
 \|\phi\|_{\infty}+2^{2j}\|{\tilde c}\|_\infty\right)
  \leq  K\,n_2 2^{2j}
\end{eqnarray*}
where $K\geq 4 \|{\tilde c}\|_\infty(\|{\tilde
c}\|_\infty+3L^4\|\phi\|_\infty^2)$. Moreover,
\begin{eqnarray*}
C^2 & =&n_2^2 \sum_{k,p} \E_{\tilde c} \left[ \left(\phi_{j,k}(U_{i_1},V_{i_1}) -
\E_{\tilde c}[ \phi_{j,k}(U,V )] \right) \left(\phi_{j,p}(U_{i_1},V_{i_1}) -
\E_{\tilde c}[
\phi_{j,p}(U,V) ] \right) \right]  \\
 & &  \times \E_{\tilde c} \left[
\left(\phi_{j,k}(U_{i_2},V_{i_2}) - \E_{\tilde c}[ \phi_{j,k}(U,V )] \right)
\left(\phi_{j,p}(U_{i_2},V_{i_2}) - \E_{\tilde c}[ \phi_{j,p}(U,V) ] \right)
\right] \\&&\\ & =& n_2^2 \sum_{k,p}\left(\E_{\tilde c} \left[
\phi_{j,k}(U_{i_1},V_{i_1})\phi_{j,p}(U_{i_1},V_{i_1})\right]-\E_{\tilde c}
\left[\phi_{j,k}(U,V)\right]\E_{\tilde c}\left[\phi_{j,p}(U,V)\right]\right)^2\\
&=&n_2^2 \sum_{k,p}\left(
\int \phi_{j,k} \phi_{j,p}\, {\tilde c} - \int \phi_{j,k}\,{\tilde c} \int \phi_{j,p}\,{\tilde c} \right)^2\\
&\leq &n_2^2 \sum_{k,p}\left(
\int \phi_{j,k} \phi_{j,p}\, {\tilde c} \right)^2+n_2^2 \left(\sum_{k}\left(\int \phi_{j,k}\,{\tilde c} \right)^2\right)^2\\
&\leq &n_2^2 \sum_{k}
\int \phi_{j,k}^2 \, {\tilde c}^2 +n_2^2 \left(\int {\tilde c}^2\right)^2\\
& \leq & \|{\tilde c}\|_\infty^2n_2^2 2^{2j} \,+ n_2^2\,\|{\tilde c}\|_2^4\leq K \,n_2^2
2^{2j},
\end{eqnarray*}
where $K\geq 2\|{\tilde c}\|_\infty^2$. Denote
$u_{\Omega_1,\Omega_2}=\E_{\tilde c}[\Omega_{\tilde c}(Z_{1},Z_{2})
\Omega_{1,{\tilde c}}(Z_{1})
 \Omega_{2,{\tilde c}}(Z_{2})]$ and for $i=1,\;2$, put
 $$
 c_i(k)=\int (\phi_{j,k}  -\E{\tilde c}
\phi_{j,k} (U,V)) \Omega_{i,{\tilde c}}\, {\tilde c}.
$$ By H\"older
Inequality, we get
\begin{eqnarray*}
u_{\Omega_1,\Omega_2} & = & \sum_k \left(\int (\phi_{j,k}  -\E_{\tilde c}
\phi_{j,k} (U,V))
\Omega_{1,{\tilde c}}\, {\tilde c}\right)\left(\int (\phi_{j,k} -\E_{\tilde c}  \phi_{j,k} (U,V))  \Omega_{2,{\tilde c}}\,{\tilde c}\right) \nonumber \\
 & \leq & \sqrt{\sum_k (c_1(k))^2 \sum_k (c_2(k))^2 } \label{FirstD}.
\end{eqnarray*}
Applying again the inequality of H\"older to $ \sum_k (c_1(k))^2$
(the same occurs for $c_2(k)$), one gets
\begin{eqnarray*}
\sum_k (c_1(k))^2 \!\!\!\!& \leq & \!\!\!\sum_k (\int (\phi_{j,k}  -\E_{\tilde c}
\phi_{j,k} (U,V))
\Omega_{1,{\tilde c}}\, {\tilde c}  \I_{[\frac{k_1}{2^{j}}, \frac{2N-1+k_1}{2^{j}}] \times [\frac{k_2}{2^{j}}, \frac{2N-1+k_2}{2^{j}}]})^2 \nonumber \\
& \leq & \sum_k \left( \int (\phi_{j,k}  -\E_{\tilde c}
\phi_{j,k} (U,V))^2 {\tilde c}  \right)  \times \nonumber \\
& & \hspace{0.5cm}\left(\int (\Omega_{1, {\tilde c}})^2\, \I_{[\frac{k_1}{2^{j}}, \frac{2L-1+k_1}{2^{j}}] \times [\frac{k_2}{2^{j}}, \frac{2L-1+k_2}{2^{j}}]} \; {\tilde c}    \right) \nonumber \\
& \leq & \|\tilde c\|_\infty  \int (\Omega_{1, {\tilde c}})^2\, {\tilde c}  \sum_k \I_{[\frac{k_1}{2^{j}}, \frac{2L-1+k_1}{2^{j}}] \times [\frac{k_2}{2^{j}}, \frac{2L-1+k_2}{2^{j}}]}\nonumber \\
\label{SecondD} \!\!\!\! & \leq & \!\!\!\! 12 \|\phi\|_{\infty}^2
L^2,
\end{eqnarray*}
since $\E_{\tilde c} (\Omega_{1, {\tilde c}}(U))^2) \leq 1$. It
follows that $D \leq K\,n_2$, where $K > 12 L^2\|\phi\|_{\infty}^2$.

\bibliographystyle{plain}

\bibliography{CEK}

\begin{table}[t!]
\begin{center}
\begin{tabular}{llccc}
\hline &   &   &   &\\
Copula \hspace{0.2cm} & True copula & $\tau=0.25$& $\tau=0.50$&
$\tau=0.75$
\\ under $H_0$&   &   &   &\\
&   &   &   &\\
\hline
&   &   &   &\\
Gumbel&Clayton&  1.00     {\it     0.0000} (0.72)& 1.00  {\it  0.0000} (0.99) &  1.00 {\it 0.0000} (1.00)\\
       &Gumbel&{\bf 0.05   {\it 0.0105 }    (0.05)}&{\bf 0.01    0.0049  (0.05)}&{\bf 0.00  {\it  0.0000}  (0.05)}\\
        &Frank&  0.31    {\it 0.0207 } (0.15)&   0.36   {\it  0.0214 } (0.40 )  & 0.24  {\it  0.0191} (0.84)\\
        &Normal&    0.25  {\it   0.0195} (0.10)&  0.25  {\it   0.0194 } (0.18)& 0.09   {\it 0.0127}    (0.61)\\
        &Student(4)&  0.99  {\it   0.0035} (0.14)& 0.89   {\it  0.0139 } (0.22)&  0.62  {\it  0.0218} (0.55)\\
        &   &   &   &\\\hline&   &   &   &\\
Clayton&Clayton&{\bf 0.05 {\it    0.0101} (0.05)}& {\bf 0.51  {\it  0.0224} (0.05)} &{\bf 1.00 {\it 0.0000}  (0.05)}\\
       &Gumbel& 0.99   {\it 0.0028} (0.86)& 1.00  {\it  0.0000} (1.00)& 1.00 {\it 0.0000} (1.00)\\
        &Frank& 0.97    {\it 0.0079} (0.56)&1.00  {\it  0.0000} (0.96)&  1.00 {\it 0.0000} (1.00)\\
        &Normal& 0.77 {\it 0.0190} (0.50)& 0.85 {\it  0.0178} (0.93)&  1.00 {\it 0.0000} (1.00)\\
        &Student(4)& 0.45  {\it  0.0223} (0.56)&0.14   {\it  0.0157} (0.95)& 1.00 {\it 0.0000} (1.00)\\
        &   &   &   &\\\hline&   &   &   &\\
Frank&Clayton&  0.99 {\it 0.0022} (0.40)& 1.00 {\it 0.0000} (0.89)  & 1.00 {\it 0.000} (0.97)\\
       &Gumbel& 0.19   {\it  0.0175} (0.33)&0.23 {\it 0.0190}  (0.63)& 0.22  {\it 0.0184} (0.82)\\
        &Frank&{\bf 0.05  {\it  0.0108 } (0.05)}&{\bf 0.01  {\it 0.0035 } (0.05)}&{\bf 0.00 {\it 0.000} (0.05)}\\
        &Normal& 0.14  {\it 0.0155}  (0.08)&0.39     {\it 0.0218} (0.20)& 0.72     {\it 0.0201} (0.41)\\
        &Student(4)& 0.95  {\it 0.0096}  (0.18)&  0.83   {\it 0.0167}   (0.08)& 0.92     {\it 0.0121} (0.06)\\
        &   &   &   &\\\hline&   &   &   &\\
Normal&Clayton& 0.97  {\it    0.0076}  (0.31)& 1.00 {\it 0.0000} (0.80) & 1.00 {\it 0.0000} (0.92)\\
       &Gumbel& 0.19  {\it    0.0176} (0.24)&  0.13 {\it 0.0151}  (0.38)& 0.01  {\it  0.0040} (0.38)\\
        &Frank& 0.23  {\it    0.0190}  (0.08)&  0.35  {\it  0.0214}  (0.20)& 0.60 {\it 0.0219}     (0.42)\\
        &Normal&{\bf  0.05   {\it   0.0099}    (0.05)}&{\bf  0.01   {\it   0.0045}    (0.05)}&{\bf 0.00 {\it 0.0000} (0.05)}\\
        &Student(4)&   0.87   {\it   0.0149}  (0.10)&  0.22  {\it   0.0185} (0.08)& 0.08 {\it 0.0120}  (0.06)\\
    &   &   &   &\\\hline&   &   &   &\\
Student(4)&Clayton&  0.71  {\it  0.0204}   (0.27)&  1.00 {\it 0.000}  (0.77)& 1.00 {\it 0.000}   (0.93)\\
       &Gumbel&  0.98   {\it 0.0056} (0.19)& 0.74 {\it  0.0331}   (0.34)&  0.29 {\it 0.0202 }(0.42)\\
        &Frank& 0.28 {\it  0.4485} (0.09)& 0.80   {\it  0.0182}   (0.27)&  0.02  {\it  0.0061 }(0.41)\\
        &Normal&  0.84 {\it 0.0166}  (0.05)&   0.20   {\it  0.0178 } (0.04)& 0.03   {\it 0.0076 }(0.04)\\
        &Student(4)&{\bf 0.03  {\it  0.0074} (0.05)}&{\bf 0.01 {\it  0.0034} (0.05)}&{\bf 0.00 {\it 0.0000}   (0.05)}\\
&   &   &   &\\\hline
\end{tabular}\end{center}
\hspace{0.5cm}
 \caption{$n_{MC}=500,n_B=20,n=2048,nn=2048$. Seed $1$. Empirical power for the test of $H_0: c=c_{\lambda_0}$ at the given level
$\alpha=10\%$ where $c_{\lambda_0}$ is specified in the first column
and the data are issue from a copula density specified in the second
column. The parameter of each copula density is chosen such that the
Kendall's tau
 is respectively $\tau=0.25,0.50,0.75$.
 }
 \label{outside}
\end{table}

\newpage

\begin{table}[h!]\begin{center}
\begin{tabular}{ccccc}
Family&parameter grid&Cardinal& $\hat \alpha$ &
Decision\\
\hline
 Gumbel&$1.05:0.1:1.95$&10&0.00  &Yes\\
Gaussian&$0.0:0.1:0.9$&10&0.04 &Yes\\
 Clayton&$0.5:0.1:1.4$&10&0.42&Yes\\
 Frank&$1.5:0.5:6.0$&10& 1.00&No\\
\hline
 Gumbel&$1.0:0.05:1.95$&20&0.00 &Yes\\
Gaussian&$0.0:0.05:0.95$&20&0.00   &Yes\\
 Clayton&$0.5:0.05:1.45$&20&0.54  &No\\
 Frank&$1.25:0.25:6.0$&20&1.00&No\\
\hline
 Gumbel &$1.45$&1&0.10&Yes\\
Gaussien&$0.48$&1&0.12&Yes\\
Clayton&$0.92$&1&0.86&No\\
Frank&$3.20$&1&1.00    &No\\
\hline
 Gumbel &$1.36\quad(1.07\%)$&1&0.02&Yes\\
Gaussien&$0.45\quad(3.27\%)$&1&0.08&Yes\\
Clayton& $0.41\;(13.15\%)$&1&0.62  &No\\
 Frank&  $2.88\quad(3.93\%)$&1&1.00     &No\\
\hline
\end{tabular}
\end{center}
\hspace{0.5cm} \caption{Empirical probability $\hat \alpha$ to
reject the fit to a fixed parametrical family given in the first
column and Decision at the prescribed level $\alpha=5\%$.  Multivariate null hypotheses (first
and second part);  $H_0 \; :\: c=c_{\hat\lambda}$, where
$\hat \lambda$ is obtained by inversion of the empirical Kendall's tau (third part);
$H_0 \; :\; c=c_{\tilde\lambda}$, where $\tilde \lambda$ is obtained by
minimizing the ASE quantity which is given into brackets (fourth part).}
 \label{test_frees}
\end{table}

\end{document}